\documentclass[11pt]{article}
\usepackage{epic,latexsym,amssymb,xcolor}
\usepackage{color}
\usepackage{tikz}
\usepackage{amsfonts,epsf,amsmath,leftidx}
\usepackage{float}
\usepackage{pgfplots}
\usepackage{mathrsfs}
\usetikzlibrary{arrows}
\usepackage{pdflscape}
\usepackage[bottom=1.5in, top=1.2in, left=1.5in, right=1.5in]{geometry}

\newtheorem{thm}{Theorem}

\newcommand{\qed}{$\Box$}

\let\oldenumerate\enumerate
\renewcommand{\enumerate}{
  \oldenumerate
  \setlength{\itemsep}{0pt}
  \setlength{\parskip}{0pt}
  \setlength{\parsep}{0pt}
}

\begin{document}

\title{The Expected Values of Hosoya Index and Merrifield-Simmons Index of Random Hexagonal Cacti}

\author{Moe Moe Oo$^{a,c}$, Nathakhun Wiroonsri$^{a, c}$, Natawat Klamsakul$^{a, c}$ \\ Thiradet Jiarasuksakun$^{b,c}$, Pawaton Kaemawichanurat$^{a, c}$
\\ \\
$^a$Department of Mathematics, Faculty of Science,\\
King Mongkut's University of Technology Thonburi, \\
Bangkok, Thailand \\
$^{b}$The Institute for the Promotion of Teaching Science and Technology (IPST),\\ Bangkok, Thailand\\
$^c$Mathematics and Statistics with Applications (MaSA) \\
\small \tt moe.oo01@mail.kmutt.ac.th, nathakhun.wir@mail.kmutt.ac.th, natawat.kla@kmutt.ac.th, \\
\small \tt thiradet@ipst.ac.th, pawaton.kae@kmutt.ac.th}

\date{}
\maketitle

\begin{abstract}
Hosoya index and Merrifield-Simmons index are two well-known topological descriptors that reflex some physical properties, such as boiling points and heat of formation, of bezenoid hydrocarbon compounds. In this paper, we establish the generating functions of the expected values of these two indices of random hexagonal cacti. This generalizes the results of Doslic and Maloy, published in Discrete Mathematics in 2010.  By applying the ideas on meromorphic functions and the growth of power series coefficients, the asymptotic behaviors of these indices on the random cacti have been established.
\end{abstract}


\section{Introduction and Motivation}
Throughout this paper, we may let $G = (V, E)$ be a graph having the vertex set $V$ and the edge set $E$. For a vertex $v \in V$, we say that $v$ is a \emph{cut vertex} of $G$ if $G - v$ has more components than $G$. A maximal connected subgraph $H$ of $G$ such that $H$ does not have a cut vertex is called a \emph{block}. Hence, if $G$ has a cut vertex, then $H$ contains a cut vertex of $G$. A block $B$ is an \emph{end block} if $B$ has exactly one cut vertex of $G$. A vertex subset $I \subseteq V$ is \emph{independent} if any pair of vertices in $I$ are not adjacent in $G$. An edge subset $M \subseteq E$ is \emph{matching} if any two edges in $M$ are not adjacent to a common vertex. A \emph{cycle} of length $k$ is denoted by $C_{k}$. A \emph{regular hexagonal cactus} is a graph that has exactly two end blocks and all the blocks are $C_{6}$, hexagons. A regular hexagonal cactus $G$ is said to be \emph{ortho} if the two cut vertices of $G$ that belong to the same non-end block are adjacent. A regular hexagonal cactus $G$ is said to be \emph{meta} if the two cut vertices of $G$ that belong to the same non-end block are a pair of vertices at distance two. Further, a regular hexagonal cactus $G$ is said to be \emph{para} if the two cut vertices of $G$ that belong to the same non-end block is a pair of vertices at distance three. For a non-negative integer $n$ and non-negative real numbers $a, b, c$ such that $a + b + c = 1$, a \emph{random hexagonal cactus} $R_{n}(a, b, c)$ with $n$ hexagons is defined as follows: when $n = 0$, $R_{0}(a, b, c)$ is the empty graph. When $n = 1$, $R_{1}(a, b, c)$ is a hexagon and, when $n = 2$, $R_{2}(a, b, c)$ is obtained from two hexagons by identifying one vertex of each. For $n \geq 3$, renaming hexagons if necessary, we label the names of hexagons of $R_{n - 1}(a, b, c)$ by $1, ..., n - 1$ consecutively along the cactus. Further, the graph $R_{n}(a, b, c)$ is obtained from $R_{n - 1}(a, b, c)$ and a hexagon $H$ by identifying a vertex of $H$ with a vertex at distance one with the probability $a$ or a vertex at distance two with the probability $b$ or the vertex at distance three with the probability $c$ from the cut vertex of the $(n - 1)^{th}$ hexagon of $R_{n - 1}(a, b, c)$.
\vskip 5 pt


\indent Husimi \cite{H} expanded Mayer and Mayer \cite{Mayer}'s book of  Statistical Mechanics by generalizing  cluster and irreducible integrals for the Theory of Condensation in the year 1952. Husimi's integrals can be represented using graphs by indication of Uhlenbeck \cite{U} in the following same year. These kinds of graphs were known as  Husimi trees which are graphs with each edge is in at most one cycle. In the fact that Husimi trees can clarify many of condensation phenomena, for example of the studies see \cite{HN,HU,R},they have been drawn a lot of attention. Cacti (Husimi trees) were first described in graph theory literature in 1973 which was 20 years since initial introduction after Harary and Palmer \cite{HP} printed their classical book on graph enumeration 
\vskip 5 pt


\indent In 1971, Hosoya introduced, in his classical paper \cite{Ho}, a graph parameter called $Z$-index which is the total number of matchings of the graphs. Hosoya found that $Z$-index relates with boiling point of the graphs representing saturated hydrocarbons. Interestingly, Gutman et. al. \cite{GVH} found further that $Z$-index also relates with the importance molecular graph descriptor called \emph{graph energy}. This has attracted much attentions of graph theorists and has resulted in many studies to evaluate this index of graphs which has been well known in Hosoya index later. For some example of studies to find Hosoya index see \cite{GZ,WG}.
\vskip 5 pt

\indent From Merrifield and Simmons \cite{MS1,MS2,MS3,MS4}'s observation, physical properties of hydrocarbonsare related to topological indices on graphs applied in chemistry which describe molecular structures such as the number of independent sets and irredundant sets. Specifically, in \cite{MS1}, the observation from Merrifield and Simmons is that the numbers of independent sets of graphs representing Alkanes which varies inversely to the boiling points and heat of formations of the compounds. From then on, Merrifield-Simmons index are recognized as the number of independent sets of graphs representing molecular structures which have been researched by several graph theorists, see \cite{DL,DM,OC} for example.
\vskip 5 pt

\indent In 2010, Doslic and Maloy \cite{DM} established the generating functions of Hosoya index and Merrifield-Simmons index of ortho-, meta, and para hexagonal cacti as follows.

\begin{thm}\cite{DM}\label{DoM1}
For a non-negative integer $n$, let $O(x), M(x)$ and $P_{x}$ be the generating functions of the number of matchings of ortho-, meta-, para-hexagonal cacti, respectively. Then 
\begin{align*}
    O(x) &= \frac{1 + 7x}{1 - 11x - 26x^{2}}\\
    M(x) &= \frac{1 + 5x}{1 - 13x + 10x^{2}}\\
    P(x) &= \frac{1 + 6x}{1 - 12x - 8x^{2}}.
\end{align*}
\end{thm}

\begin{thm}\cite{DM}\label{DoM2}
For a non-negative integer $n$, let $\overline{O}(x), \overline{M}(x)$ and $\overline{P}_{x}$ be the generating functions of the number of independent sets of ortho-, meta-, para-hexagonal cacti, respectively. Then 
\begin{align*}
    \overline{O}(x) &= \frac{2 + 2x}{1 - 8x - 25x^{2}}\\
    \overline{M}(x) &= \frac{2 - 6x}{1 - 12x + 11x^{2}}\\
    \overline{P}(x) &= \frac{2 - 2x}{1 - 10x - 7x^{2}}.
\end{align*}
\end{thm}

\indent For related works on finding Hosoya Index and Merrrifild-Simmons Index in random chain, Huang et.al.\cite{HKD} described exact formulas for the expected values of these indices of a random polyphenylene chain including $n$ octagons. Chen et.al.\cite{A} established explicit expressions
for the expected value of Merrifield-Simmons index of a random phenylene chain $PH_n,p$ and a random hexagonal
chain $HS_np$ by using the method of the generating functions. The corresponding entropy constants are computed and the maximum and minimum values are obtained in both random systems. Very recently in 2022, Sun et.al. \cite{SG} found the recurrence relations of the expected values of Hosoya Index and Merrrifild-Simmons Index of a random cyclooctylene chains containing $n$ octagons. By solving these recurrence relations, the formula of expected values of these indices were established. 
\vskip 5 pt

\indent In this paper, we establish generating function of the expected values of Hosoya index and Merrifield-Simmons index of random hexagonal cacti. Some special cases of our results prove Theorems \ref{DoM1} and \ref{DoM2}.

\section{Main Results}
In this section, we present our theorems related to the expected values of Hosoya Index and Merrifiled-Simmonds Index of random hexagonal cacti. Our first theorem establishes generating function of expected values of Hosoya index of random hexagonal cacti
\vskip 5 pt

\begin{thm}\label{theorem1}
For a non-negative integer $n$ and non-negative real numbers $a, b, c$ such that $a + b + c = 1$, we let $R_{n}(a, b, c)$ be a random hexagonal chain cactus with $n$ hexagons. Further, we let $E(m_{n}(a, b, c))$ be the expected value of the number of matchings of $R_{n}(a, b, c)$ and $M_{a, b, c}(x)$ be the generating function of $E(m_{n}(a, b, c))$. Then 
\begin{align*}
M_{a, b, c}(x) = \sum_{n = 0}^{\infty}E(m_{n}(a, b, c))x^{n} = \frac{1 + 10x - 3ax -5bx - 4cx}{1 - 8x -3ax - 5bx - 4cx - 26ax^{2} + 10bx^{2} - 8cx^{2}}.
\end{align*}
\end{thm}

\noindent By letting $a = 1, b = c = 0$ and $a = c = 0, b = 1$ and $a = b = 0, c = 1$ the graph $R_{n}(a, b, c)$ becomes ortho-, meta- and para-hexagonal cacti, respectively. Thus, we have the following equations:
\begin{align*}
    M_{1, 0, 0}(x) &= O(x) = \frac{1 + 7x}{1 - 11x - 26x^{2}}\\
    M_{0, 1, 0}(x) &= M(x) = \frac{1 + 5x}{1 - 13x + 10x^{2}}\\
    M_{0, 0, 1}(x) &= P(x) = \frac{1 + 6x}{1 - 12x - 8x^{2}}.
\end{align*}

\noindent Hence, Theorem \ref{theorem1} generalizes Theorem \ref{DoM1}.

\indent Further, by considering the principal part of the expansion of the generating function in Theorem \ref{theorem1} around the singularity,  we have the asymptotic behavior of $E(m_{n}(a, b, c))$ as follows:

\begin{thm}\label{asymp-t}
If $E(m_{n}(a, b, c))$ is the expected value of the number of matchings of $R_{n}(a, b, c)$. Then,
\begin{align*}
   E(m_{n}(a, b, c))\approx \frac{\frac{1}{(\frac{1}{4})^{n + 1}}(\sigma_{2} - 3a\sigma_{1} - 5b\sigma_{1} - 4c\sigma_{1} + 9a^{2} + 25b^{2} + 16c^{2} + 10\sigma_{1} - 80)}{(-\frac{3a + 5b + 4c - \sigma_{1} + 8}{13a - 5b + 4c})^{n + 1}(52a - 20b + 16c)\sigma_{1}}. 
\end{align*}
\noindent where $\sigma_{1} = \sqrt{9a^{2} + 30ab + 24ac + 152a + 25b^{2} + 40bc + 40b + 16c^{2} + 96c + 64}$ and $\sigma_{2} = 46a - 30b + 8c + 30ab + 24ac + 40bc$.
\end{thm}

\begin{thm}\label{theorem2}
For a non-negative integer $n$ and non-negative real numbers $a, b, c$ such that $a + b + c = 1$, we let $R_{n}(a, b, c)$ be a random hexagonal chain cactus with $n$ hexagons. Further, we let $E(i_{n}(a, b, c))$ be the expected value of the number of independent sets of $R_{n}(a, b, c)$ and $I_{a, b, c}(x)$ be the generating function of $E(i_{n}(a, b, c))$. Then 
\begin{align*}
I_{a, b, c}(x) = \sum_{n = 0}^{\infty}E(i_{n}(a, b, c))x^{n} = \frac{1 + 13x - 3ax - 7bx - 5cx + 25ax^{2} - 11bx^{2} + 7cx^{2}}{1 - 5x - 3ax - 7bx - 5cx - 25ax^{2} + 11bx^{2} - 7cx^{2}}.
\end{align*}
\end{thm}

\noindent By letting $a = 1, b = c = 0$ and $a = c = 0, b = 1$ and $a = b = 0, c = 1$ the graph $R_{n}(a, b, c)$ becomes ortho-, meta- and para-hexagonal cacti, respectively. Thus, we have the following equations:
\begin{align*}
    I_{1, 0, 0}(x) &= \overline{O}(x) = \frac{1 + 10x + 25x^{2}}{1 - 8x - 25x^{2}}\\
    I_{0, 1, 0}(x) &= \overline{M}(x) = \frac{1 + 6x -11x^{2}}{1 - 12x + 11x^{2}}\\
    I_{0, 0, 1}(x) &= \overline{P}(x) = \frac{1 + 8x + 7x^{2}}{1 - 10x - 7x^{2}}.
\end{align*}

\noindent Hence, Theorem \ref{theorem2} completes Theorem \ref{DoM2}.

\indent Further, by considering the principal part of the expansion of the generating function in Theorem \ref{theorem2} around the singularity, we have the asymptotic behavior of $E(i_{n}(a, b, c))$ as follows:

\begin{thm}\label{asymp-i}
If $E(i_{n}(a, b, c))$ is the expected value of the number of independent sets of $R_{n}(a, b, c)$. Then,
\begin{align*}
E(i_{n}(a, b, c))\approx \frac{\frac{1}{0.5^{n + 1}}(- 3a\sigma_{1} - 7b3a\sigma_{1} - 5c3a\sigma_{1} + 42ab + 30ac + 70bc + 4\sigma_{1} + \sigma_{3})}{(-\frac{3a + 7b + 5c - \sigma_{1} + 5}{\sigma_{2}})^{n + 1}\sigma_{1}\sigma_{2}}. 
\end{align*}
\noindent where $\sigma_{1} = \sqrt{9a^{2} + 42ab + 30ac + 130a + 49b^{2} + 70bc + 26b + 25c^{2} + 78c + 25}$ and $\sigma_{2} = 25a - 11b + 7c$ and $\sigma_{3} = 9a^{2} + 49b^{2} + 25c^{2} - 20 + 53a - 15b + 19c $.
\end{thm}

\section{Proofs}
\begin{figure}[h!]
\centering

\tikzset{every picture/.style={line width=0.75pt}} 

\begin{tikzpicture}[x=0.75pt,y=0.75pt,yscale=-1,xscale=1]

\draw   (155,112.13) -- (142.5,133.78) -- (117.5,133.78) -- (105,112.13) -- (117.5,90.47) -- (142.5,90.47) -- cycle ;
\draw   (10.5,112.13) .. controls (10.5,95.9) and (31.65,82.75) .. (57.75,82.75) .. controls (83.85,82.75) and (105,95.9) .. (105,112.13) .. controls (105,128.35) and (83.85,141.5) .. (57.75,141.5) .. controls (31.65,141.5) and (10.5,128.35) .. (10.5,112.13) -- cycle ;
\draw   (118.17,40.48) -- (139.65,53.27) -- (139.32,78.27) -- (117.5,90.47) -- (96.02,77.68) -- (96.36,52.69) -- cycle ;
\draw   (324,109.13) -- (311.5,130.78) -- (286.5,130.78) -- (274,109.13) -- (286.5,87.47) -- (311.5,87.47) -- cycle ;
\draw   (179.5,109.13) .. controls (179.5,92.9) and (200.65,79.75) .. (226.75,79.75) .. controls (252.85,79.75) and (274,92.9) .. (274,109.13) .. controls (274,125.35) and (252.85,138.5) .. (226.75,138.5) .. controls (200.65,138.5) and (179.5,125.35) .. (179.5,109.13) -- cycle ;
\draw   (495,109.13) -- (482.5,130.78) -- (457.5,130.78) -- (445,109.13) -- (457.5,87.47) -- (482.5,87.47) -- cycle ;
\draw   (350.5,109.13) .. controls (350.5,92.9) and (371.65,79.75) .. (397.75,79.75) .. controls (423.85,79.75) and (445,92.9) .. (445,109.13) .. controls (445,125.35) and (423.85,138.5) .. (397.75,138.5) .. controls (371.65,138.5) and (350.5,125.35) .. (350.5,109.13) -- cycle ;
\draw   (545,109.13) -- (532.5,130.78) -- (507.5,130.78) -- (495,109.13) -- (507.5,87.47) -- (532.5,87.47) -- cycle ;
\draw   (349,65.82) -- (336.5,87.47) -- (311.5,87.47) -- (299,65.82) -- (311.5,44.17) -- (336.5,44.17) -- cycle ;


\draw (120,109.38) node   [align=left] {\begin{minipage}[lt]{10.88pt}\setlength\topsep{0pt}
{\fontfamily{ptm}\selectfont \textit{\small{$H_{n - 1}$}}}
\end{minipage}};

\draw (290,109.38) node   [align=left] {\begin{minipage}[lt]{10.88pt}\setlength\topsep{0pt}
{\fontfamily{ptm}\selectfont \textit{\small{$H_{n - 1}$}}}
\end{minipage}};

\draw (460,109.38) node   [align=left] {\begin{minipage}[lt]{10.88pt}\setlength\topsep{0pt}
{\fontfamily{ptm}\selectfont \textit{\small{$H_{n - 1}$}}}
\end{minipage}};

\draw (520,109.38) node   [align=left] {\begin{minipage}[lt]{10.88pt}\setlength\topsep{0pt}
{\fontfamily{ptm}\selectfont \textit{\small{$H_{n}$}}}
\end{minipage}};

\draw (120,70) node   [align=left] {\begin{minipage}[lt]{10.88pt}\setlength\topsep{0pt}
{\fontfamily{ptm}\selectfont \textit{\small{$H_{n}$}}}
\end{minipage}};

\draw (320,70) node   [align=left] {\begin{minipage}[lt]{10.88pt}\setlength\topsep{0pt}
{\fontfamily{ptm}\selectfont \textit{\small{$H_{n}$}}}
\end{minipage}};

\end{tikzpicture}

\vskip -0.25 cm
\caption{\footnotesize{The graphs $R_{n}(a, b, c)$}}
\label{Rn}
\end{figure}
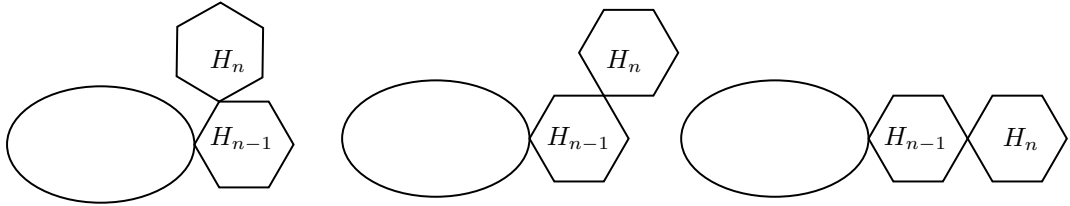

In this section, we prove Theorems \ref{theorem1} and \ref{theorem2}. Firstly, we recall that, for non-negative real numbers $a, b, c$ such that $a + b + c = 1$, the graph $R_{0}(a, b, c)$ is empty, the graph $R_{1}(a, b, c)$ is a hexagon and the graph $R_{2}(a, b, c)$ is obtained from two hexagons by identifying one vertex of each. For $n \geq 3$, the graph $R_{n}(a, b, c)$ is obtained from $R_{n - 1}(a, b, c)$ and a hexagon $H$ by identifying a vertex of $H$ with a vertex at distance one with the probability $a$ or a vertex at distance two with the probability $b$ or the vertex at distance three with the probability $c$ from the cut vertex of the $(n - 1)^{th}$ hexagon of $R_{n - 1}(a, b, c)$. Figure \ref{Rn} shows examples of $R_{n}(a, b, c)$.

\vskip 20 pt

\noindent In the figure, from left to right, a vertex of the hexagon $H_{n}$ is identified with a vertex of the hexagon $H_{n - 1}$ with the probabilities $a, b$ and $c$, respectively.
\vskip 5 pt

\begin{figure}[h!]
\centering

\tikzset{every picture/.style={line width=0.75pt}} 

\begin{tikzpicture}[x=0.75pt,y=0.75pt,yscale=-1,xscale=1]

\draw   (155,112.13) -- (142.5,133.78) -- (117.5,133.78) -- (105,112.13) -- (117.5,90.47) -- (142.5,90.47) -- cycle ;
\draw   (10.5,112.13) .. controls (10.5,95.9) and (31.65,82.75) .. (57.75,82.75) .. controls (83.85,82.75) and (105,95.9) .. (105,112.13) .. controls (105,128.35) and (83.85,141.5) .. (57.75,141.5) .. controls (31.65,141.5) and (10.5,128.35) .. (10.5,112.13) -- cycle ;
\draw   (324,109.13) -- (311.5,130.78) -- (286.5,130.78) -- (274,109.13) -- (286.5,87.47) -- (311.5,87.47) -- cycle ;
\draw   (179.5,109.13) .. controls (179.5,92.9) and (200.65,79.75) .. (226.75,79.75) .. controls (252.85,79.75) and (274,92.9) .. (274,109.13) .. controls (274,125.35) and (252.85,138.5) .. (226.75,138.5) .. controls (200.65,138.5) and (179.5,125.35) .. (179.5,109.13) -- cycle ;
\draw   (495,109.13) -- (482.5,130.78) -- (457.5,130.78) -- (445,109.13) -- (457.5,87.47) -- (482.5,87.47) -- cycle ;
\draw   (350.5,109.13) .. controls (350.5,92.9) and (371.65,79.75) .. (397.75,79.75) .. controls (423.85,79.75) and (445,92.9) .. (445,109.13) .. controls (445,125.35) and (423.85,138.5) .. (397.75,138.5) .. controls (371.65,138.5) and (350.5,125.35) .. (350.5,109.13) -- cycle ;
\draw    (95.37,77.67) -- (117.05,90.12) ;
\draw    (95.37,77.67) -- (95.33,52.67) ;
\draw    (116.95,40.13) -- (95.33,52.67) ;
\draw    (138.63,52.58) -- (116.95,40.13) ;

\draw    (299.01,65.98) -- (312.24,87.2) ;
\draw    (299.01,65.98) -- (310.77,43.92) ;
\draw    (335.76,43.07) -- (310.77,43.92) ;
\draw    (348.98,64.29) -- (335.76,43.07) ;

\draw    (507.57,87.94) -- (494.98,109.54) ;
\draw    (507.57,87.94) -- (532.57,88.06) ;
\draw    (544.98,109.76) -- (532.57,88.06) ;
\draw    (532.38,131.36) -- (544.98,109.76) ;

\draw (130,109.38) node   [align=left] {\begin{minipage}[lt]{10.88pt}\setlength\topsep{0pt}
{\fontfamily{ptm}\selectfont \textit{\small{$H_{n}$}}}
\end{minipage}};

\draw (300,109.38) node   [align=left] {\begin{minipage}[lt]{10.88pt}\setlength\topsep{0pt}
{\fontfamily{ptm}\selectfont \textit{\small{$H_{n}$}}}
\end{minipage}};

\draw (470,109.38) node   [align=left] {\begin{minipage}[lt]{10.88pt}\setlength\topsep{0pt}
{\fontfamily{ptm}\selectfont \textit{\small{$H_{n}$}}}
\end{minipage}};

\end{tikzpicture}

\vskip -0.25 cm
\caption{\footnotesize{The graphs $R'_{n}(a, b, c)$}}
\label{Rprimen}
\end{figure}
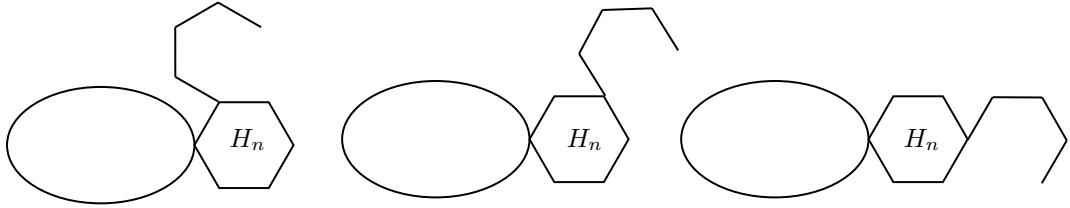

\indent We further define three auxiliary graphs as follows. The graph $R'_{n}(a, b, c)$ is obtained from $R_{n}(a, b, c)$ and a path of length four $P_{5} = x_{1}x_{2}x_{3}x_{4}x_{5}$ by identifying $x_{1}$ of $P_{5}$ to a vertex at distance one with the probability $a$, a vertex at distance two with the probability $b$ and the vertex at distance three with the probability $c$ from the cut vertex of the $n^{th}$ hexagon of $R_{n}(a, b, c)$. Figure \ref{Rprimen} illustrates examples of $R'_{n}(a, b, c)$. From left to right, $x_{1}$ of the path $x_{1}...x_{5}$ is identified with a vertex of the hexagon $H_{n}$ with the probabilities $a, b$ and $c$, respectively.
\vskip 15 pt

\vskip 20 pt

\begin{figure}[h!]
\centering

\tikzset{every picture/.style={line width=0.75pt}} 

\begin{tikzpicture}[x=0.75pt,y=0.75pt,yscale=-1,xscale=1]

\draw   (155,112.13) -- (142.5,133.78) -- (117.5,133.78) -- (105,112.13) -- (117.5,90.47) -- (142.5,90.47) -- cycle ;
\draw   (10.5,112.13) .. controls (10.5,95.9) and (31.65,82.75) .. (57.75,82.75) .. controls (83.85,82.75) and (105,95.9) .. (105,112.13) .. controls (105,128.35) and (83.85,141.5) .. (57.75,141.5) .. controls (31.65,141.5) and (10.5,128.35) .. (10.5,112.13) -- cycle ;
\draw   (324,109.13) -- (311.5,130.78) -- (286.5,130.78) -- (274,109.13) -- (286.5,87.47) -- (311.5,87.47) -- cycle ;
\draw   (179.5,109.13) .. controls (179.5,92.9) and (200.65,79.75) .. (226.75,79.75) .. controls (252.85,79.75) and (274,92.9) .. (274,109.13) .. controls (274,125.35) and (252.85,138.5) .. (226.75,138.5) .. controls (200.65,138.5) and (179.5,125.35) .. (179.5,109.13) -- cycle ;
\draw   (495,109.13) -- (482.5,130.78) -- (457.5,130.78) -- (445,109.13) -- (457.5,87.47) -- (482.5,87.47) -- cycle ;
\draw   (350.5,109.13) .. controls (350.5,92.9) and (371.65,79.75) .. (397.75,79.75) .. controls (423.85,79.75) and (445,92.9) .. (445,109.13) .. controls (445,125.35) and (423.85,138.5) .. (397.75,138.5) .. controls (371.65,138.5) and (350.5,125.35) .. (350.5,109.13) -- cycle ;
\draw    (139.36,52.78) -- (117.8,40.13) ;
\draw    (139.36,52.78) -- (139.18,77.78) ;
\draw    (117.44,90.12) -- (139.18,77.78) ;
\draw    (95.88,77.47) -- (117.44,90.12) ;

\draw    (348.61,64.56) -- (335.34,43.37) ;
\draw    (348.61,64.56) -- (336.88,86.64) ;
\draw    (311.9,87.53) -- (336.88,86.64) ;
\draw    (298.64,66.34) -- (311.9,87.53) ;

\draw    (533.04,130.48) -- (545.64,108.89) ;
\draw    (533.04,130.48) -- (508.04,130.36) ;
\draw    (495.64,108.65) -- (508.04,130.36) ;
\draw    (508.25,87.06) -- (495.64,108.65) ;

\draw (130,109.38) node   [align=left] {\begin{minipage}[lt]{10.88pt}\setlength\topsep{0pt}
{\fontfamily{ptm}\selectfont \textit{\small{$H_{n}$}}}
\end{minipage}};

\draw (300,109.38) node   [align=left] {\begin{minipage}[lt]{10.88pt}\setlength\topsep{0pt}
{\fontfamily{ptm}\selectfont \textit{\small{$H_{n}$}}}
\end{minipage}};

\draw (470,109.38) node   [align=left] {\begin{minipage}[lt]{10.88pt}\setlength\topsep{0pt}
{\fontfamily{ptm}\selectfont \textit{\small{$H_{n}$}}}
\end{minipage}};

\end{tikzpicture}

\vskip -0.25 cm
\caption{\footnotesize{The graphs $\tilde{R}_{n}(a, b, c)$}}
\label{Rtilden}
\end{figure}
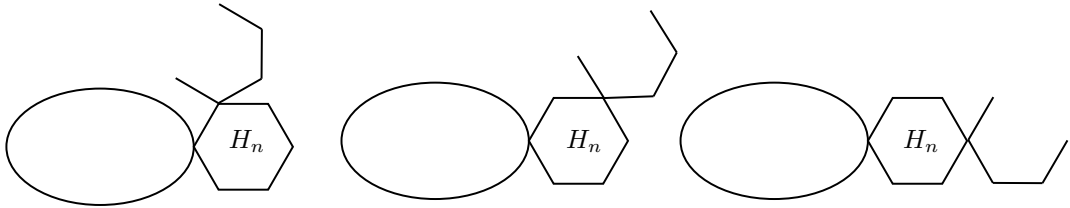

\noindent The graph $\tilde{R}_{n}(a, b, c)$ is obtained from $R_{n}(a, b, c)$ and a path of length four $P_{5} = x_{1}x_{2}x_{3}x_{4}x_{5}$ by identifying $x_{2}$ of $P_{5}$ to a vertex at distance one with the probability $a$, a vertex at distance two with the probability $b$ and the vertex at distance three with the probability $c$ from the cut vertex of the $n^{th}$ hexagon of $R_{n}(a, b, c)$. Figure \ref{Rtilden} shows examples of $\tilde{R}_{n}(a, b, c)$. From left to right, $x_{2}$ of the path $x_{1}...x_{5}$ is identified with a vertex of the hexagon $H_{n}$ with the probabilities $a, b$ and $c$, respectively.
\vskip 15 pt

\vskip 20 pt

\begin{figure}[h!]
\centering

\tikzset{every picture/.style={line width=0.75pt}} 

\begin{tikzpicture}[x=0.75pt,y=0.75pt,yscale=-1,xscale=1]

\draw   (155,112.13) -- (142.5,133.78) -- (117.5,133.78) -- (105,112.13) -- (117.5,90.47) -- (142.5,90.47) -- cycle ;
\draw   (10.5,112.13) .. controls (10.5,95.9) and (31.65,82.75) .. (57.75,82.75) .. controls (83.85,82.75) and (105,95.9) .. (105,112.13) .. controls (105,128.35) and (83.85,141.5) .. (57.75,141.5) .. controls (31.65,141.5) and (10.5,128.35) .. (10.5,112.13) -- cycle ;
\draw   (324,109.13) -- (311.5,130.78) -- (286.5,130.78) -- (274,109.13) -- (286.5,87.47) -- (311.5,87.47) -- cycle ;
\draw   (179.5,109.13) .. controls (179.5,92.9) and (200.65,79.75) .. (226.75,79.75) .. controls (252.85,79.75) and (274,92.9) .. (274,109.13) .. controls (274,125.35) and (252.85,138.5) .. (226.75,138.5) .. controls (200.65,138.5) and (179.5,125.35) .. (179.5,109.13) -- cycle ;
\draw   (495,109.13) -- (482.5,130.78) -- (457.5,130.78) -- (445,109.13) -- (457.5,87.47) -- (482.5,87.47) -- cycle ;
\draw   (350.5,109.13) .. controls (350.5,92.9) and (371.65,79.75) .. (397.75,79.75) .. controls (423.85,79.75) and (445,92.9) .. (445,109.13) .. controls (445,125.35) and (423.85,138.5) .. (397.75,138.5) .. controls (371.65,138.5) and (350.5,125.35) .. (350.5,109.13) -- cycle ;
\draw    (139.48,77.2) -- (139.48,52.2) ;
\draw    (139.48,77.2) -- (117.83,89.7) ;
\draw    (96.18,77.2) -- (117.83,89.7) ;
\draw    (96.18,52.2) -- (96.18,77.2) ;

\draw    (336.65,86.79) -- (348.64,64.85) ;
\draw    (336.65,86.79) -- (311.65,87.37) ;
\draw    (298.65,66.02) -- (311.65,87.37) ;
\draw    (310.64,44.08) -- (298.65,66.02) ;

\draw    (508.54,130.64) -- (533.54,130.21) ;
\draw    (508.54,130.64) -- (495.67,109.21) ;
\draw    (507.79,87.35) -- (495.67,109.21) ;
\draw    (532.79,86.92) -- (507.79,87.35) ;

\draw (130,109.38) node   [align=left] {\begin{minipage}[lt]{10.88pt}\setlength\topsep{0pt}
{\fontfamily{ptm}\selectfont \textit{\small{$H_{n}$}}}
\end{minipage}};

\draw (300,109.38) node   [align=left] {\begin{minipage}[lt]{10.88pt}\setlength\topsep{0pt}
{\fontfamily{ptm}\selectfont \textit{\small{$H_{n}$}}}
\end{minipage}};

\draw (470,109.38) node   [align=left] {\begin{minipage}[lt]{10.88pt}\setlength\topsep{0pt}
{\fontfamily{ptm}\selectfont \textit{\small{$H_{n}$}}}
\end{minipage}};

\end{tikzpicture}

\vskip -0.25 cm
\caption{\footnotesize{The graphs $\hat{R}_{n}(a, b, c)$}}
\label{Rhatn}
\end{figure}
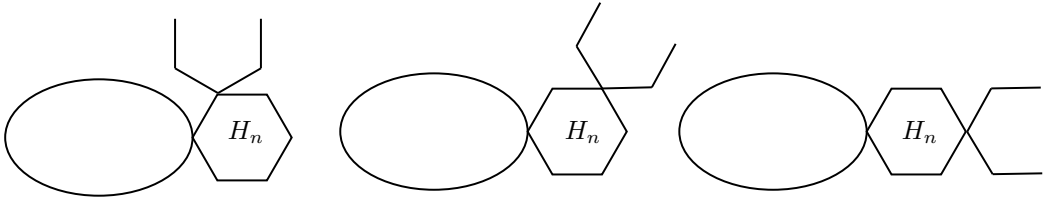

\noindent The graph $\hat{R}_{n}(a, b, c)$ is obtained from $R_{n}(a, b, c)$ and a path of length four $P_{5} = x_{1}x_{2}x_{3}x_{4}x_{5}$ by identifying $x_{3}$ of $P_{5}$ to a vertex at distance one with the probability $a$, a vertex at distance two with the probability $b$ and the vertex at distance three with the probability $c$ from the cut vertex of the $n^{th}$ hexagon of $R_{n}(a, b, c)$. Figure \ref{Rhatn} illustrates examples of $\hat{R}_{n}(a, b, c)$. From left to right, $x_{3}$ of the path $x_{1}...x_{5}$ is identified with a vertex of the hexagon $H_{n}$ with the probabilities $a, b$ and $c$, respectively.
\vskip 15 pt

\vskip 20 pt

\subsection{Hosoya Index of Random Hexagoanl Chain Cacti}
In this subsection, we prove Theorem \ref{theorem1}. Recall that 
\vskip 5 pt

\indent $E(m_{n}(a, b, c))$ is the expected value of the number of matchings of $R_{n}(a, b, c)$
\vskip 5 pt

\indent $M_{a, b, c}(x)$ is the generating function of $E(m_{n}(a, b, c))$. 
\vskip 5 pt

\noindent When there is no danger of confusion, we let $E(m_{n})$ to denote $E(m_{n}(a, b, c))$ and let $M(x)$ to denote $M_{a, b, c}(x)$. Thus,
\vskip 5 pt

$M(x) = \sum^{\infty}_{n = 0}E(m_{n})x^{n}$.

\noindent Similarly, we let $E(m'_{n}), E(\tilde{m}_{n})$ and $E(\hat{m}_{n})$ be the expected values of the number of matchings of $R'_{n}(a, b, c), \tilde{R}_{n}(a, b, c)$ and $\hat{R}_{n}(a, b, c)$, respectively. Further, we let $M'(x), \tilde{M}(x)$ and $\hat{M}(x)$ be the generating functions of $E(m'_{n}), E(\tilde{m}_{n})$ and $E(\hat{m}_{n})$, respectively. Therefore,
\vskip 5 pt

\indent $M'(x) = \sum^{\infty}_{n = 0}E(m'_{n})x^{n}$
\vskip 5 pt

\indent $\tilde{M}(x) = \sum^{\infty}_{n = 0}E(\tilde{m}_{n})x^{n}$ 
\vskip 5 pt

\indent $\hat{M}(x) = \sum^{\infty}_{n = 0}E(\hat{m}_{n})x^{n}$.
\vskip 5 pt

\indent Now, to prove Theorem \ref{theorem1}, we will prove the following equations.
\begin{align}
    (1 - 8x)M(x)          &= 1 +10x + 10ax^{2}M'(x) + 10bx^{2}\tilde{M}(x) + 10cx^{2}\hat{M}(x)\label{eqm1}\\
    (1 - 3ax)M'(x)        &= 3 + 5M(x) + 3bx\tilde{M}(x) + 3cx\hat{M}(x)\label{eqm2}\\
    (1 - 5bx)\tilde{M}(x) &= 5 + 3M(x) + 5axM(x) + 3cx\hat{M}(x)\label{eqm3}\\
    (1 - 4cx)\hat{M}(x)   &= 4 + 4M(x) + 4axM(x) + 4bx\tilde{M}(x)\label{eqm4}
\end{align}

\noindent \emph{Proof of Equation} (\ref{eqm1}). We name the vertices of the $n^{th}$ hexagon $H_{n}$ of $R_{n}(a, b, c)$ by $h_{1}, h_{2}, ..., h_{6}$ clockwise with $h_{1}$ is the cut vertex containing in $H_{n}$ (and in the $(n - 1)^{th}$ hexagon $H_{n - 1}$). Similarly, we name the vertices of $H_{n - 1}$ by $k_{1}, k_{1}, ..., k_{6}$ clockwise with $k_{1}$ is the cut vertex containing in $H_{n - 1}$. Because $R_{n}(a, b, c)$ is a chain, $k_{1} \neq h_{1}$. We distinguish three cases depending on the distance between $h_{1}$ and $k_{1}$ on $H_{n - 1}$.
\vskip 5 pt

\noindent \textbf{Case 1:} $h_{1}$ is adjacent to $k_{1}$.\\
\indent Thus, $h_{1}$ is either $k_{2}$ or $k_{6}$. Without loss of generality, we let $h_{1} = k_{2}$. This case occurs with the probability $a$. We count the number of matchings due to the following subcases.
\vskip 5 pt

\noindent \textbf{Subscase 1.1:} The matchings that contain $h_{1}h_{2}$.\\
\indent For any matching that contains $h_{1}h_{2}$, it cannot contain $k_{2}k_{3}, k_{1}k_{2}, h_{2}h_{3}$ and $h_{1}h_{6}$. Thus, the number of matchings in this subcase equals the number of matchings of the union of $R'_{n - 2}(a, b, c)$ and the path $P_{4} = h_{3}h_{4}h_{5}h_{6}$. Clearly, $P_{4}$ has $5$ matchings. So, the expected value of the number of matchings in this subcase is $5E(m'_{n - 2})$.
\vskip 5 pt

\noindent \textbf{Subscase 1.2:} The matchings that contain $h_{1}h_{6}$.\\
\indent This subcase is symmetric to Subcase 1.1. Thus,  the expected value of the number of matchings in this subcase is $5E(m'_{n - 2})$.
\vskip 5 pt

\noindent \textbf{Subscase 1.3:} The matchings that contain neither $h_{1}h_{2}$ nor $h_{1}h_{6}$.\\
\indent  The number of matchings in this subcase equals the number of matchings of the union of $R_{n - 1}(a, b, c)$ and the path $P_{5} = h_{2}h_{3}h_{4}h_{5}h_{6}$. Clearly, $P_{5}$ has $8$ matchings. So, the expected value of the number of matchings in this subcase is $8E(m_{n - 1})$.
\vskip 5 pt

\indent From the three subcases, we have that there are $10aE(m'_{n - 2}) + 8aE(m_{n - 1})$.
\vskip 5 pt

\noindent \textbf{Case 2:} $h_{1}$ is at distance two from $k_{1}$.\\
\indent Thus, $h_{1}$ is either $k_{3}$ or $k_{5}$. Without loss of generality, we let $h_{1} = k_{3}$. This case occurs with the probability $b$. We count the number of matchings due to the following subcases.
\vskip 5 pt

\noindent \textbf{Subscase 2.1:} The matchings that contain $h_{1}h_{2}$.\\
\indent For any matching that contains $h_{1}h_{2}$, it cannot contain $k_{2}k_{3}, k_{3}k_{4}, h_{2}h_{3}$ and $h_{1}h_{6}$. Thus,  the number of matchings in this subcase equals the number of matchings of the union of $\tilde{R}_{n - 2}(a, b, c)$ and the path $P_{4} = h_{3}h_{4}h_{5}h_{6}$. Clearly, $P_{4}$ has $5$ matchings. So, the expected value of the number of matchings in this subcase is $5E(\tilde{m}_{n - 2})$.
\vskip 5 pt

\noindent \textbf{Subscase 2.2:} The matchings that contain $h_{1}h_{6}$.\\
\indent This subcase is symmetric to Subcase 2.1. Thus,  the expected value of the number of matchings in this subcase is $5E(\tilde{m}_{n - 2})$.
\vskip 5 pt

\noindent \textbf{Subscase 2.3:} The matchings that contain neither $h_{1}h_{2}$ nor $h_{1}h_{6}$.\\
\indent  The number of matchings in this subcase equals the number of matchings of the union of $R_{n - 1}(a, b, c)$ and the path $P_{5} = h_{2}h_{3}h_{4}h_{5}h_{6}$. Clearly, $P_{5}$ has $8$ matchings. So, the expected value of the number of matching in this subcase is $8E(m_{n - 1})$.
\vskip 5 pt

\indent From the three subcases, we have that there are $10bE(\tilde{m}_{n - 2}) + 8bE(m_{n - 1})$.
\vskip 5 pt

\noindent \textbf{Case 3:} $h_{1}$ is at distance three from $k_{1}$.\\
\indent Thus, $h_{1}$ is $k_{4}$. This case occurs with the probability $c$. We count the number of matchings due to the following subcases.
\vskip 5 pt

\noindent \textbf{Subscase 3.1:} The matchings that contain $h_{1}h_{2}$.\\
\indent For any matching that contains $h_{1}h_{2}$, it cannot contain $k_{4}k_{5}, k_{3}k_{4}, h_{2}h_{3}$ and $h_{1}h_{6}$.  the number of matchings in this subcase equals the number of matchings of the union of $\hat{R}_{n - 2}(a, b, c)$ and a path $P_{4} = h_{3}h_{4}h_{5}h_{6}$. Clearly, $P_{4}$ has $5$ matchings. So, the expected value of the number of matchings in this subcase is $5E(\hat{m}_{n - 2})$.
\vskip 5 pt

\noindent \textbf{Subscase 3.2:} The matchings that contain $h_{1}h_{6}$.\\
\indent This subcase is symmetric to Subcase 3.1. Thus, the expected value of the number of matchings in this subcase is $5E(\hat{m}_{n - 2})$.
\vskip 5 pt

\noindent \textbf{Subscase 3.3:} The matchings that contain neither $h_{1}h_{2}$ nor $h_{1}h_{6}$.\\
\indent  The number of matchings in this subcase equals the number of matchings of the union of $R_{n - 1}(a, b, c)$ and the path $P_{5} = h_{2}h_{3}h_{4}h_{5}h_{6}$. Clearly, $P_{5}$ has $8$ matchings. So, the expected value of the number of matchings in this subcase is $8E(m_{n - 1})$.
\vskip 5 pt

\indent From the three subcases, we have that there are $10cE(\hat{m}_{n - 2}) + 8cE(m_{n - 1})$. Therefore, from Cases $1, 2$ and $3$, we have that
\begin{align*}
E(m_{n}) &= 10aE(m'_{n - 2}) + 8aE(m_{n - 1}) + 10bE(\tilde{m}_{n - 2}) + 8bE(m_{n - 1}) + 10cE(\hat{m}_{n - 2}) + 8cE(m_{n - 1})\notag\\
         &= 8(a + b + c)E(m_{n - 1}) + 10aE(m'_{n - 2}) + 10bE(\tilde{m}_{n - 2}) + 10cE(\hat{m}_{n - 2})\notag\\
         &= 8E(m_{n - 1}) + 10aE(m'_{n - 2}) + 10bE(\tilde{m}_{n - 2}) + 10cE(\hat{m}_{n - 2}).\notag
\end{align*}

\noindent For $n \geq 2$, we multiply $x^{n}$ throughout the above equation and sum over all $n$. We have that
\begin{align*}
\sum^{\infty}_{n = 2}E(m_{n})x^{n} = \sum^{\infty}_{n = 2}8E(m_{n - 1})x^{n} + 10a\sum^{\infty}_{n = 2}E(m'_{n - 2})x^{n} + 10b\sum^{\infty}_{n = 2}E(\tilde{m}_{n - 2})x^{n} + 10c\sum^{\infty}_{n = 2}E(\hat{m}_{n - 2})x^{n}\notag
\end{align*}
\noindent which implies that
\begin{align*}
M(x) - E(m_{0}) - E(m_{1})x = 8x(M(x) - E(m_{0})) + 10ax^{2}M'(x) + 10bx^{2}\tilde{M}(x) + 10cx^{2}\hat(M)(x).\notag
\end{align*}
It can be checked that $E(m_{0}) = 1$ and $E(m_{1}) = 18$. Hence, 
\begin{align*}
M(x) = 1 + 10x + 8xM(x) + 10ax^{2}M'(x) + 10bx^{2}\tilde{M}(x) + 10cx^{2}\hat(M)(x).
\end{align*}
This proves Equation (\ref{eqm1}).
\qed

\noindent \emph{Proof of Equation} (\ref{eqm2}). Recall that $R'_{n}(a, b, c)$ is obtained from $R_{n}(a, b, c)$ and a path $P_{5} = x_{1}x_{2}x_{3}x_{4}x_{5}$ by identifying $x_{1}$ with a vertex of the $n^{th}$ hexagon of $R_{n}(a, b, c)$. We name the vertices of the $n^{th}$ hexagon $H_{n}$ of $R'_{n}(a, b, c)$ by $h_{1}, h_{2}, ..., h_{6}$ clockwise with $h_{1}$ is the cut vertex containing in $H_{n}$ (and in the $(n - 1)^{th}$ hexagon $H_{n - 1}$). By the definition of $R'_{n}(a, b, c)$, we have that $x_{1} \neq h_{1}$. We distinguish three cases due to the distance between $x_{1}$ and $h_{1}$. 
\vskip 5 pt

\noindent \textbf{Case 1:} $x_{1}$ is adjacent to $h_{1}$.\\
\indent Thus, $x_{1}$ is either $h_{2}$ or $h_{6}$. Without loss of generality, we let $x_{1} = h_{2}$. This case occurs with the probability $a$. We count the number of matchings due to the following subcases.
\vskip 5 pt

\noindent \textbf{Subscase 1.1:} The matchings that contain $x_{1}x_{2}$.\\
\indent For any matching that contains $x_{1}x_{2}$, it cannot contain $x_{2}x_{3}, h_{1}h_{2}$ and $h_{2}h_{3}$. Thus, the number of matchings in this subcase equals the number of matchings of the union of $R'_{n -1}(a, b, c)$ and the path $x_{3}x_{4}x_{5}$. Clearly, $x_{3}x_{4}x_{5}$ has $3$ matchings. So, the expected value of the number of matchings in this subcase is $3E(m'_{n - 1})$.
\vskip 5 pt

\noindent \textbf{Subscase 1.2:} The matchings that do not contain $x_{1}x_{2}$.\\
\indent  The number of matchings in this subcase equals the number of matchings of the union of $R_{n}(a, b, c)$ and the path $x_{2}x_{3}x_{4}x_{5}$. Clearly, $x_{2}x_{3}x_{4}x_{5}$ has $5$ matchings. So, the expected value of the number of matchings in this subcase is $5E(m_{n})$.
\vskip 5 pt

\indent From the two subcases, we have that there are $3aE(m'_{n - 1}) + 5aE(m_{n})$ from Case 1.
\vskip 5 pt

\noindent \textbf{Case 2:} $x_{1}$ is at distance two from $h_{1}$.\\
\indent Thus, $x_{1}$ is either $h_{3}$ or $h_{5}$. Without loss of generality, we let $x_{1} = h_{3}$. This case occurs with the probability $b$. We count the number of matchings due to the following subcases.
\vskip 5 pt

\noindent \textbf{Subscase 2.1:} The matchings that contain $x_{1}x_{2}$.\\
\indent For any matching that contains $x_{1}x_{2}$, it cannot contain $x_{2}x_{3}, h_{3}h_{4}$ and $h_{2}h_{3}$. Thus,  the number of matchings in this subcase equals the number of matchings of the union of $\tilde{R}_{n - 1}(a, b, c)$ and the path $x_{3}x_{4}x_{5}$. Clearly, $x_{3}x_{4}x_{5}$ has $3$ matchings. So, the expected value of the number of matchings in this subcase is $3E(\tilde{m}_{n - 1})$.
\vskip 5 pt

\noindent \textbf{Subscase 2.2:} The matchings that do not contain $x_{1}x_{2}$.\\
\indent  The number of matchings in this subcase equals the number of matchings of the union of $R_{n}(a, b, c)$ and the path $x_{2}x_{3}x_{4}x_{5}$. Clearly, $x_{2}x_{3}x_{4}x_{5}$ has $5$ matchings. So, the expected value of the number of matchings in this subcase is $5E(m_{n})$.
\vskip 5 pt

\indent From the two subcases, we have that there are $3bE(\tilde{m}_{n - 1}) + 5bE(m_{n})$ from Case 2.
\vskip 5 pt

\noindent \textbf{Case 3:} $x_{1}$ is at distance three from $h_{1}$.\\
\indent Thus, $x_{1}$ is $h_{4}$. This case occurs with the probability $c$. We count the number of matchings due to the following subcases.
\vskip 5 pt

\noindent \textbf{Subscase 3.1:} The matchings that contain $x_{1}x_{2}$.\\
\indent For any matching that contains $x_{1}x_{2}$, it cannot contain $x_{2}x_{3}, h_{3}h_{4}$ and $h_{4}h_{5}$. Thus,  the number of matchings in this subcase equals the number of matchings of the union of $\hat{R}_{n - 1}(a, b, c)$ and the path $x_{3}x_{4}x_{5}$. Clearly, $x_{3}x_{4}x_{5}$ has $3$ matchings. So, the expected value of the number of matching in this subcase is $3E(\hat{m}_{n - 1})$.
\vskip 5 pt

\noindent \textbf{Subscase 3.2:} The matchings that do not contain $x_{1}x_{2}$.\\
\indent  The number of matchings in this subcase equals the number of matchings of the union of $R_{n}(a, b, c)$ and the path $x_{2}x_{3}x_{4}x_{5}$. Clearly, $x_{2}x_{3}x_{4}x_{5}$ has $5$ matchings. So, the expected value of the number of matchings in this subcase is $5E(m_{n})$.
\vskip 5 pt

\indent From the two subcases, we have that there are $3cE(\hat{m}_{n - 2}) + 5cE(m_{n})$ from Case 3. Therefore, from Cases $1, 2$ and $3$, we have that
\begin{align*}
E(m'_{n}) &= 3aE(m'_{n - 1}) + 5aE(m_{n}) + 3bE(\tilde{m}_{n - 1}) + 5bE(m_{n}) + 3cE(\hat{m}_{n - 1}) + 5cE(m_{n})\notag\\
          &= 5(a + b + c)E(m_{n}) + 3aE(m'_{n - 1}) + 3bE(\tilde{m}_{n - 1}) + 3cE(\hat{m}_{n - 1})\notag\\
          &= 5E(m_{n}) + 3aE(m'_{n - 1}) + 3bE(\tilde{m}_{n - 1}) + 3cE(\hat{m}_{n - 1}).\notag
\end{align*}

\noindent For $n \geq 1$, we multiply $x^{n}$ throughout the above equation and sum over all $n$. We have that
\begin{align*}
\sum^{\infty}_{n = 1}E(m'_{n})x^{n} = \sum^{\infty}_{n = 1}5E(m_{n})x^{n} + 3a\sum^{\infty}_{n = 1}E(m'_{n - 1})x^{n} + 3b\sum^{\infty}_{n = 1}E(\tilde{m}_{n - 1})x^{n} + 3c\sum^{\infty}_{n = 1}E(\hat{m}_{n - 1})x^{n}\notag
\end{align*}
\noindent which implies that
\begin{align*}
M'(x) - E(m'_{0}) = 5(M(x) - E(m_{0})) + 3axM'(x) + 3bx\tilde{M}(x) + 3cx\hat{M}(x).\notag
\end{align*}
It can be checked that $E(m'_{0}) = 8$ and $E(m_{0}) = 1$. Hence, 
\begin{align*}
M'(x) = 3 + 5M(x) + 3axM'(x) + 3bx\tilde{M}(x) + 3cx\hat{M}(x).
\end{align*}
This proves Equation (\ref{eqm2}).
\qed

\noindent \emph{Proof of Equation} (\ref{eqm3}). Recall that $\tilde{R}_{n}(a, b, c)$ is obtained from $R_{n}(a, b, c)$ and a path $P_{5} = x_{1}x_{2}x_{3}x_{4}x_{5}$ by identifying $x_{2}$ with a vertex of the $n^{th}$ hexagon of $R_{n}(a, b, c)$. We name the vertices of the $n^{th}$ hexagon $H_{n}$ of $\tilde{R}_{n}(a, b, c)$ by $h_{1}, h_{2}, ..., h_{6}$ clockwise with $h_{1}$ is the cut vertex containing in $H_{n}$ (and in the $(n - 1)^{th}$ hexagon $H_{n - 1}$). By the definition of $\tilde{R}_{n}(a, b, c)$, we have that $x_{2} \neq h_{1}$. We distinguish three cases due to the distance between $x_{2}$ and $h_{1}$. 
\vskip 5 pt

\noindent \textbf{Case 1:} $x_{2}$ is adjacent to $h_{1}$.\\
\indent Thus, $x_{2}$ is either $h_{2}$ or $h_{6}$. Without loss of generality, we let $x_{2} = h_{2}$. This case occurs with the probability $a$. We count the number of matchings due to the following subcases.
\vskip 5 pt

\noindent \textbf{Subscase 1.1:} The matchings that contain $x_{1}x_{2}$.\\
\indent For any matching that contains $x_{1}x_{2}$, it cannot contain $x_{2}x_{3}, h_{1}h_{2}$ and $h_{2}h_{3}$. Thus,  the number of matchings in this subcase equals the number of matchings of the union of $R'_{n - 1}(a, b, c)$ and the path $P_{3} = x_{3}x_{4}x_{5}$. Clearly, $P_{3}$ has $3$ matchings. So, the expected value of the number of matchings in this subcase is $3E(m'_{n - 1})$.
\vskip 5 pt

\noindent \textbf{Subscase 1.2:} The matchings that contain $x_{2}x_{3}$.\\
\indent For any matching that contains $x_{2}x_{3}$, it cannot contain $x_{1}x_{2}, x_{3}x_{4}, h_{1}h_{2}$ and $h_{2}h_{3}$. Thus,  the number of matchings in this subcase equals the number of matchings of the union of $R'_{n - 1}(a, b, c)$, the vertex $x_{1}$ and the path $P_{2} = x_{4}x_{5}$. Clearly, $P_{2}$ has $2$ matchings. So, the expected value of the number of matchings in this subcase is $2E(m'_{n - 1})$.
\vskip 5 pt

\noindent \textbf{Subscase 1.3:} The matchings that contain neither $x_{1}x_{2}$ nor $x_{2}x_{3}$.\\
\indent  The number of matchings in this subcase equals the number of matchings of the union of $R_{n}(a, b, c)$ and the path $P_{3} = x_{3}x_{4}x_{5}$. Clearly, $P_{3}$ has $3$ matchings. So, the expected value of the number of matchings in this subcase is $3E(m_{n})$.
\vskip 5 pt

\indent From the three subcases, we have that there are $5aE(m'_{n - 1}) + 3aE(m_{n})$.
\vskip 5 pt

\noindent \textbf{Case 2:} $x_{2}$ is at distance two from $h_{1}$.\\
\indent Thus, $x_{2}$ is either $h_{3}$ or $h_{5}$. Without loss of generality, we let $x_{2} = h_{3}$. This case occurs with the probability $b$. We count the number of matchings due to the following subcases.
\vskip 5 pt

\noindent \textbf{Subscase 2.1:} The matchings that contain $x_{1}x_{2}$.\\
\indent For any matching that contains $x_{1}x_{2}$, it cannot contain $x_{2}x_{3}, h_{2}h_{3}$ and $h_{3}h_{4}$. Thus,  the number of matchings in this subcase equals the number of matchings of the union of $\tilde{R}_{n - 1}(a, b, c)$ and the path $P_{3} = x_{3}x_{4}x_{5}$. Clearly, $P_{3}$ has $3$ matchings. So, the expected value of the number of matchings in this subcase is $3E(\tilde{m}_{n - 1})$.
\vskip 5 pt

\noindent \textbf{Subscase 2.2:} The matchings that contain $x_{2}x_{3}$.\\
\indent For any matching that contains $x_{2}x_{3}$, it cannot contain $x_{1}x_{2}, x_{3}x_{4}, h_{3}h_{4}$ and $h_{2}h_{3}$. Thus,  the number of matchings in this subcase equals the number of matchings of the union of $\tilde{R}_{n - 1}(a, b, c)$, the vertex $x_{1}$ and the path $P_{2}= x_{4}x_{5}$. Clearly, $P_{2}$ has $2$ matchings. So, the expected value of the number of matchings in this subcase is $2E(\tilde{m}_{n - 1})$.
\vskip 5 pt

\noindent \textbf{Subscase 2.3:} The matchings that contain neither $x_{1}x_{2}$ nor $x_{2}x_{3}$.\\
\indent  The number of matchings in this subcase equals the number of matchings of the union of $R_{n}(a, b, c)$ and the path $P_{3} = x_{3}x_{4}x_{5}$. Clearly, $P_{3}$ has $3$ matchings. So, the expected value of the number of matchings in this subcase is $3E(m_{n})$.
\vskip 5 pt

\indent From the three subcases, we have that there are $5bE(\tilde{m}_{n - 1}) + 3bE(m_{n})$.
\vskip 5 pt

\noindent \textbf{Case 3:} $x_{2}$ is at distance three from $h_{1}$.\\
\indent Thus, $x_{2}$ is $h_{4}$. This case occurs with the probability $c$. We count the number of matchings due to the following subcases.
\vskip 5 pt

\noindent \textbf{Subscase 3.1:} The matchings that contain $x_{1}x_{2}$.\\
\indent For any matching that contains $x_{1}x_{2}$, it cannot contain $x_{2}x_{3}, h_{3}h_{4}$ and $h_{4}h_{5}$. Thus, the number of matchings in this subcase equals the number of matchings of the union of $\hat{R}_{n - 1}(a, b, c)$ and the path $P_{3}= x_{3}x_{4}x_{5}$. Clearly, $P_{3}$ has $3$ matchings. So, the expected value of the number of matchings in this subcase is $3E(\hat{m}_{n - 1})$.
\vskip 5 pt

\noindent \textbf{Subscase 3.2:} The matchings that contain $x_{2}x_{3}$.\\
\indent For any matching that contains $x_{2}x_{3}$, it cannot contain $x_{1}x_{2}, x_{3}x_{4}, h_{3}h_{4}$ and $h_{4}h_{5}$. Thus,  the number of matchings in this subcase equals the number of matchings of the union of $\hat{R}_{n - 1}(a, b, c)$, the vertex $x_{1}$ and the path $P_{2} =x_{4}x_{5}$. Clearly, $P_{2}$ has $2$ matchings. So, the expected value of the number of matchings in this subcase is $2E(\hat{m}_{n - 1})$.
\vskip 5 pt

\noindent \textbf{Subscase 3.3:} The matchings that contain neither $x_{1}x_{2}$ nor $x_{2}x_{3}$.\\
\indent  The number of matchings in this subcase equals the number of matchings of the union of $R_{n}(a, b, c)$ and the path $P_{3} = x_{3}x_{4}x_{5}$. Clearly, $P_{3}$ has $3$ matchings. So, the expected value of the number of matchings in this subcase is $3E(m_{n})$.
\vskip 5 pt

\indent From the three subcases, we have that there are $5cE(\hat{m}_{n - 1}) + 3cE(m_{n})$. Therefore, from Cases $1, 2$ and $3$, we have that

\begin{align*}
E(\tilde{m}_{n}) &=  5aE(m'_{n - 1}) + 3aE(m_{n})+ 5bE(\tilde{m}_{n - 1}) + 3bE(m_{n}) + 5cE(\hat{m}_{n - 1}) + 3cE(m_{n }))\notag\\
         &= 3(a + b + c)E(m_{n}) + 5aE(m'_{n - 1}) + 5bE(\tilde{m}_{n - 1}) + 5cE(\hat{m}_{n - 1}) \notag\\
         &= 3E(m_{n}) + 5aE(m'_{n - 1}) + 5bE(\tilde{m}_{n - 1}) + 5cE(\hat{m}_{n - 1}).\notag
\end{align*}

\noindent For $n \geq 1$, we multiply $x^{n}$ throughout the above equation and sum over all $n$. We have that
\begin{align*}
\sum^{\infty}_{n = 1}E(\tilde{m}_{n})x^{n} = \sum^{\infty}_{n = 1}3E(m_{n})x^{n} + 5a\sum^{\infty}_{n = 1}E(m'_{n - 1})x^{n} + 5b\sum^{\infty}_{n = 1}E(\tilde{m}_{n - 1})x^{n} + 5c\sum^{\infty}_{n = 1}E(\hat{m}_{n - 1})x^{n}\notag
\end{align*}

\noindent which implies that
\begin{align*}
\Tilde{M}(x) - E(\tilde{m}_{0}) = 3(M(x) - E(m_{0})) + 5axM'(x) + 5bx\tilde{M}(x) + 5cx\hat{M}(x).\notag
\end{align*}
It can be checked that $E(\tilde{m}_{0}) = 8$ and $E(m_{0}) = 1$. Hence, 
\begin{align*}
\Tilde{M}(x) = 5 + 3M(x) + 5axM'(x) + 5bx\tilde{M}(x) + 5cx\hat{M}(x).\\
\end{align*}
         
\noindent This proves Equation (\ref{eqm3}).
\qed
\vskip 5 pt

\noindent \emph{Proof of Equation} (\ref{eqm4}). Recall that $\hat{R}_{n}(a, b, c)$ is obtained from $R_{n}(a, b, c)$ and a path $P_{5} = x_{1}x_{2}x_{3}x_{4}x_{5}$ by identifying ${x_3}$ with a vertex of the $n^{th}$ hexagon of $R_{n}(a, b, c)$. We name the vertices of the $n^{th}$ hexagon $H_{n}$ of $\hat{R}_{n}(a, b, c)$ by $h_{1}, h_{2}, ..., h_{6}$ clockwise with $h_{1}$ is the cut vertex containing in $H_{n}$ (and in the $(n - 1)^{th}$ hexagon $H_{n - 1}$). By the definition of $\hat{R}_{n}(a, b, c)$, we have that $x_{3} \neq h_{1}$. We distinguish three cases due to the distance between $x_{3}$ and $h_{1}$. 
\vskip 5 pt

\noindent \textbf{Case 1:} $x_{3}$ is adjacent to $h_{1}$.\\
\indent Thus, $x_{3}$ is either $h_{2}$ or $h_{6}$. Without loss of generality, we let $x_{3} = h_{2}$. This case occurs with the probability $a$. We count the number of matchings due to the following subcases.
\vskip 5 pt

\noindent \textbf{Subscase 1.1:} The matchings that contain $x_{2}x_{3}$.\\
\indent For any matching that contains $x_{2}x_{3}$, it cannot contain $x_{1}x_{2}, x_{3}x_{4}, h_{1}h_{2}$ and $h_{2}h_{3}$. Thus, the number of matchings in this subcase equals the number of matchings of the union of $R'_{n - 1}(a, b, c)$, the vertex $x_{1}$ and the path $P_{2} = x_{4}x_{5}$. Clearly, $ P_{2}$ has $2$ matchings. So, the expected value of the number of matchings in this subcase is $2E(m'_{n - 1})$.
\vskip 5 pt

\noindent \textbf{Subscase 1.2:} The matchings that contain $x_{3}x_{4}$.\\
\indent This subcase is symmetric to Subcase 1.1. Thus,  the expected value of the number of matchings in this subcase is $2E(m'_{n - 1})$.
\vskip 5 pt

\noindent \textbf{Subscase 1.2:} The matchings that contain $x_{3}x_{4}$.\\
\indent For any matching that contains $x_{3}x_{4}$, it cannot contain $x_{2}x_{3}, x_{4}x_{5}, h_{1}h_{2}$ and $h_{2}h_{3}$. Thus,  the number of matchings in this subcase equals the number of matchings of the union of $R'_{n - 1}(a, b, c)$, the vertex $x_{5}$ and the path $x_{1}x_{2}$. Clearly, $x_{1}x_{2}$ has $2$ matchings. So, the expected value of the number of matchings in this subcase is $2E(m'_{n - 1})$.
\vskip 5 pt

\noindent \textbf{Subscase 1.3:} The matchings that contain neither $x_{2}x_{3}$ nor $x_{3}x_{4}$.\\
\indent  The number of matchings in this subcase equals the number of matchings of the union of $R_{n}(a, b, c)$ , the paths $ x_{1}x_{2}$ and $ x_{4}x_{5}$  . Clearly, $ x_{1}x_{2}$ and $ x_{4}x_{5}$ has $4$ matchings. So, the expected value of the number of matchings in this subcase is $4E(m_{n})$.
\vskip 5 pt

\indent From the three subcases, we have that there are $4aE(m'_{n - 1}) + 4aE(m_{n})$.
\vskip 5 pt

\noindent \textbf{Case 2:} $x_{3}$ is at distance two from $h_{1}$.\\
\indent Thus, $x_{3}$ is either $h_{3}$ or $h_{5}$. Without loss of generality, we let $x_{3} = h_{3}$. This case occurs with the probability $b$. We count the number of matchings due to the following subcases.
\vskip 5 pt

\noindent \textbf{Subscase 2.1:} The matchings that contain $x_{2}x_{3}$.\\
\indent For any matching that contains $x_{2}x_{3}$, it cannot contain $x_{1}x_{2},
x_{3}x_{4},h_{2}h_{3}$ and $h_{3}h_{4}$. Thus,  the number of matchings in this subcase equals the number of matchings of the union of $\tilde{R}_{n - 1}(a, b, c)$, the vertex $x_{1}$ and the path $x_{4}x_{5}$. Clearly, $x_{4}x_{5}$ has $2$ matchings. So, the expected value of the number of matchings in this subcase is $2E(\tilde{m}_{n - 1})$.
\vskip 5 pt

\noindent \textbf{Subscase 2.2:} The matchings that contain $x_{3}x_{4}$.\\
\indent This subcase is symmetric to Subcase 2.1. Thus,  the expected value of the number of matchings in this subcase is $2E(\tilde{m}_{n - 1})$.
\vskip 5 pt

\noindent \textbf{Subscase 2.2:} The matchings that contain $x_{3}x_{4}$.\\
\indent For any matching that contains $x_{3}x_{4}$, it cannot contain $x_{2}x_{3}, x_{4}x_{5}, h_{2}h_{3}$ and $h_{3}h_{4}$. Thus,  the number of matchings in this subcase equals the number of matchings of the union of $\tilde{R}_{n - 1}(a, b, c)$, the vertex $x_{5}$ and the path $x_{1}x_{2}$. Clearly, $x_{1}x_{2}$ has $2$ matchings. So, the expected value of the number of matchings in this subcase is $2E(\tilde{m}_{n - 1})$.
\vskip 5 pt

\noindent \textbf{Subscase 2.3:} The matchings that contain neither $x_{2}x_{3}$ nor $x_{3}x_{4}$.\\
\indent  The number of matchings in this subcase equals the number of matchings of the union of $R_{n}(a, b, c)$ , the paths $ x_{1}x_{2}$ and $ x_{4}x_{5}$  . Clearly, $ x_{1}x_{2}$ and $ x_{4}x_{5}$ has $4$ matchings. So, the expected value of the number of matchings in this subcase is $4E(m_{n})$.
\vskip 5 pt

\indent From the three subcases, we have that there are $4bE(\tilde{m}_{n - 1}) + 4bE(m_{n})$.
\vskip 5 pt

\noindent \textbf{Case 3:} $x_{3}$ is at distance three from $h_{1}$.\\
\indent Thus, $x_{3}$ is $h_{4}$. This case occurs with the probability $c$. We count the number of matchings due to the following subcases.
\vskip 5 pt

\noindent \textbf{Subscase 3.1:} The matchings that contain $x_{2}x_{3}$.\\
\indent For any matching that contains $x_{2}x_{3}$, it cannot contain $x_{1}x_{2}, x_{3}x_{4},h_{3}h_{4}$ and $h_{4}h_{5}$. Thus,  the number of matchings in this subcase equals the number of matchings of the union of $\hat{R}_{n - 1}(a, b, c)$, the vertex $x_{1}$ and the path $x_{4}x_{5}$. Clearly, $x_{4}x_{5}$ has $2$ matchings. So, the expected value of the number of matchings in this subcase is $2E(\hat{m}_{n - 1})$.
\vskip 5 pt

\indent This subcase is symmetric to Subcase 3.1. Thus,  the expected value of the number of matchings in this subcase is $2E(\hat{m}_{n - 1})$.
\vskip 5 pt

\noindent \textbf{Subscase 3.2:} The matchings that contain $x_{3}x_{4}$.\\
\indent For any matching that contains $x_{3}x_{4}$, it cannot contain $x_{2}x_{3}, x_{4}x_{5}, h_{3}h_{4}$ and $h_{4}h_{5}$. Thus,  the number of matchings in this subcase equals the number of matchings of the union of $\hat{R}_{n - 1}(a, b, c)$, the vertex $x_{5}$ and the path $x_{1}x_{2}$. Clearly, $x_{1}x_{2}$ has $2$ matchings. So, the expected value of the number of matchings in this subcase is $2E(\hat{m}_{n - 1})$.
\vskip 5 pt

\noindent \textbf{Subscase 3.3:} The matchings that contain neither $x_{2}x_{3}$ nor $x_{3}x_{4}$.\\
\indent  The number of matchings in this subcase equals the number of matchings of the union of $R_{n}(a, b, c)$ , the paths $ x_{1}x_{2}$ and $ x_{4}x_{5}$  . Clearly, $ x_{1}x_{2}$ and $ x_{4}x_{5}$ has $4$ matchings. So, the expected value of the number of matchings in this subcase is $4E(m_{n})$.
\vskip 5 pt

\indent From the three subcases, we have that there are $4cE(\hat{m}_{n - 2}) + 4cE(m_{n })$. Therefore, from Cases $1, 2$ and $3$, we have that

\begin{align*}
E(\hat{m}_{n}) &=  4aE(m'_{n - 1}) + 4aE(m_{n})+ 4bE(\tilde{m}_{n - 1}) + 4bE(m_{n}) + 4cE(\hat{m}_{n - 1}) + 4cE(m_{n }))\notag\\
         &= 4(a + b + c)E(m_{n}) + 4aE(m'_{n - 1}) + 4bE(\tilde{m}_{n - 1}) + 4cE(\hat{m}_{n - 1}) \notag\\
         &= 4E(m_{n}) + 4aE(m'_{n - 1}) + 4bE(\tilde{m}_{n - 1}) + 4cE(\hat{m}_{n - 1}).\notag
\end{align*}

\noindent For $n \geq 1$, we multiply $x^{n}$ throughout the above equation and sum over all $n$. We have that
\begin{align*}
\sum^{\infty}_{n = 1}E(\hat{m}_{n})x^{n} = \sum^{\infty}_{n = 1}4E(m_{n})x^{n} + 4a\sum^{\infty}_{n = 1}E(m'_{n - 1})x^{n} + 4b\sum^{\infty}_{n = 1}E(\tilde{m}_{n - 1})x^{n} + 4c\sum^{\infty}_{n = 1}E(\hat{m}_{n - 1})x^{n}\notag
\end{align*}

\noindent which implies that
\begin{align*}
\hat{M}(x) - E(\hat{m}_{0}) = 4(M(x) - E(m_{0})) + 4axM'(x) + 4bx\tilde{M}(x) + 4cx\hat{M}(x).\notag
\end{align*}
It can be checked that $E(\hat{m}_{0}) = 8$ and $E(m_{0}) = 1$. Hence, 
\begin{align*}
\hat{M}(x) = 4 + 4M(x) + 4axM'(x) + 4bx\tilde{M}(x) + 4cx\hat{M}(x).\\
\end{align*}

\noindent This proves Equation (\ref{eqm4}).
\qed

\indent By Equations (\ref{eqm1}), (\ref{eqm2}), (\ref{eqm3}) and (\ref{eqm4}), it can be solved that 
\begin{align*}
M(x) = \frac{1 + 10x - 3ax -5bx - 4cx}{1 - 8x -3ax - 5bx - 4cx - 26ax^{2} + 10bx^{2} - 8cx^{2}}.
\end{align*}
This proves Theorem \ref{theorem1}.

\subsection{Merrifield-Simmons Index of Random Hexagoanl Chain Cacti}
In this subsection, we prove Theorem \ref{theorem2}. Recall that 
\vskip 5 pt

\indent $E(i_{n}(a, b, c))$ is the expected value of the number of independent sets of $R_{n}(a, b, c)$
\vskip 5 pt

\indent $I_{a, b, c}(x)$ is the generating function of $E(i_{n}(a, b, c))$. 
\vskip 5 pt

\noindent When there is no danger of confusion, we let $E(i_{n})$ to denote $E(i_{n}(a, b, c))$ and let $I(x)$ to denote $I_{a, b, c}(x)$. Thus,
\vskip 5 pt

$I(x) = \sum^{\infty}_{n = 0}E(i_{n})x^{n}$.

\noindent Similarly, we let $E(i'_{n}), E(\tilde{i}_{n})$ and $E(\hat{i}_{n})$ be the expected values of the number of independent sets of $R'_{n}(a, b, c), \tilde{R}_{n}(a, b, c)$ and $\hat{R}_{n}(a, b, c)$, respectively. Further, we let $I'(x), \tilde{I}(x)$ and $\hat{I}(x)$ be the generating functions of $E(i'_{n}), E(\tilde{i}_{n})$ and $E(\hat{i}_{n})$, respectively. Therefore,
\vskip 5 pt

\indent $I'(x) = \sum^{\infty}_{n = 0}E(i'_{n})x^{n}$
\vskip 5 pt

\indent $\tilde{I}(x) = \sum^{\infty}_{n = 0}E(\tilde{i}_{n})x^{n}$ 
\vskip 5 pt

\indent $\hat{I}(x) = \sum^{\infty}_{n = 0}E(\hat{i}_{n})x^{n}$.
\vskip 5 pt

\indent Now, to prove Theorem \ref{theorem2}, we will prove the following equations.
\begin{align}
    (1 - 5x)I(x)          &= 1 +13x + 8ax^{2}I'(x) + 8bx^{2}\tilde{I}(x) + 8cx^{2}\hat{I}(x)\label{eqm5}\\
    (1 - 3ax)I'(x)        &= 8 + 5I(x) + 3bx\tilde{I}(x) + 3cx\hat{I}(x)\label{eqm6}\\
    (1 - 7bx)\tilde{I}(x) &= 10 + 3I(x) + 7axI'(x) + 7cx\hat{I}(x)\label{eqm7}\\
    (1 - 5cx)\hat{I}(x)   &= 9 + 4I(x) + 5axI'(x) + 5bx\tilde{I}(x)\label{eqm8}
\end{align}

\noindent \emph{Proof of Equation} (\ref{eqm5}). We name the vertices of the $n^{th}$ hexagon $H_{n}$ of $R_{n}(a, b, c)$ by $h_{1}, h_{2}, ..., h_{6}$ clockwise with $h_{1}$ is the cut vertex containing in $H_{n}$ (and in the $(n - 1)^{th}$ hexagon $H_{n - 1}$). Similarly, we name the vertices of $H_{n - 1}$ by $k_{1}, k_{1}, ..., k_{6}$ clockwise with $k_{1}$ is the cut vertex containing in $H_{n - 1}$. Because $R_{n}(a, b, c)$ is a chain, $k_{1} \neq h_{1}$. We distinguish three cases depending on the distance between $h_{1}$ and $k_{1}$ on $H_{n - 1}$.
\vskip 5 pt

\noindent \textbf{Case 1:} $h_{1}$ is adjacent to $k_{1}$.\\
\indent Thus, $h_{1}$ is either $k_{2}$ or $k_{6}$. Without loss of generality, we let $h_{1} = k_{2}$. This case occurs with the probability $a$. We count the number of independent sets due to the following subcases.
\vskip 5 pt

\noindent \textbf{Subscase 1.1:} The independent sets  that contain $h_{2}$ and do not contain $h_{6}$.\\
\indent For any independent set that contains $h_{2}$ and does not contain $h_{6}$, it cannot contain $h_{1}$ and $h_{3}$. Thus, the number of independent sets in this subcase equals the number of independent sets of the union of $R'_{n - 2}(a, b, c)$ and the path $P_{2} = h_{4}h_{5}$. Clearly, $P_{2}$ has $3$ independent sets. So, the expected value of the number of independent sets in this subcase is $3E(i'_{n - 2})$.
\vskip 5 pt

\noindent \textbf{Subscase 1.2:} The independent sets  that contain $h_{6}$ and do not contain $h_{2}$.\\
\indent For any independent set that contains $h_{6}$ and does not contain $h_{2}$, it cannot contain $h_{1}$ and $h_{5}$. Thus, the number of independent sets in this subcase equals the number of independent sets of the union of $R'_{n - 2}(a, b, c)$ and the path $P_{2} = h_{3}h_{4}$. Clearly, $P_{2}$ has $3$ independent sets. So, the expected value of the number of independent sets in this subcase is $3E(i'_{n - 2})$.

\vskip 5 pt

\noindent \textbf{Subscase 1.3:} The independent sets that contain both $h_{2}$ and  $h_{6}$.\\
\indent For any independent set that contains both $h_{2}$ and $h_{6}$, it cannot contain $h_{1}$,$h_{3}$ and $h_{5}$. Thus, the number of independent sets in this subcase equals the number of independent sets of the union of $R'_{n - 2}(a, b, c)$ and the vertex $h_{4}$. Clearly, the vertex $h_{4}$ has $2$ independent sets. So, the expected value of the number of indenpendent sets in this subcase is $2E(i'_{n - 2})$.
\vskip 5 pt

\noindent \textbf{Subscase 1.4:} The independent sets that contain  neither $h_{2}$ nor  $h_{6}$.\\
\indent The number of independent sets in this subcase equals the number of independent sets of the union of $R_{n - 1}(a, b, c)$ and the path $P_{3} = h_{3}h_{4}h_{5}$. Clearly, $P_{3}$ has $5$ independent sets. So, the expected value of the number of independent sets in this subcase is $5E(i_{n - 1})$.
\vskip 5 pt

\indent From the four subcases, we have that there are $8aE(i'_{n - 2}) + 5aE(i_{n - 1})$.
\vskip 5 pt

\noindent \textbf{Case 2:} $h_{1}$ is at distance two from $k_{1}$.\\
\indent Thus, $h_{1}$ is either $k_{3}$ or $k_{5}$. Without loss of generality, we let $h_{1} = k_{3}$. This case occurs with the probability $b$. We count the number of independent sets due to the following subcases.
\vskip 5 pt

\noindent \textbf{Subscase 2.1:} The independent sets that contain $h_{2}$ and do not contain $h_{6}$.\\
\indent For any independent sets that contain $h_{2}$ and does not contain $h_{6}$, it cannot contain $h_{1}$ and $h_{3}$. Thus, the number of independent sets in this subcase equals the number of independent sets of the union of $\tilde{R}_{n - 2}(a, b, c)$ and the path $P_{2} = h_{4}h_{5}$. Clearly, $P_{2}$ has $3$ independent sets. So, the expected value of the number of independent sets in this subcase is $3E(\tilde{i}_{n - 2})$.
\vskip 5 pt

\noindent \textbf{Subscase 2.2:} The independent sets that contain $h_{6}$ and do not contain $h_{2}$.\\
\indent For any independent sets that contain $h_{6}$ and does not contain $h_{2}$, it cannot contain $h_{1}$ and $h_{5}$. Thus, the number of independent sets in this subcase equals the number of independent sets of the union of $\tilde{R}_{n - 2}(a, b, c)$ and the path $P_{2} = h_{3}h_{4}$. Clearly, $P_{2}$ has $3$ independent sets. So, the expected value of the number of independent sets in this subcase is $3E(\tilde{i}_{n - 2})$.
\vskip 5 pt

\noindent \textbf{Subscase 2.3:} The independent sets that contain both $h_{2}$ and  $h_{6}$.\\
\indent For any independent set that contains both $h_{2}$ and $h_{6}$, it cannot contain $h_{1}$,$h_{3}$ and $h_{5}$. Thus, the number of independent sets in this subcase equals the number of independent sets of the union of $\tilde{R}_{n - 2}(a, b, c)$ and the vertex $h_{4}$. Clearly, the vertex $h_{4}$ has $2$ independent sets. So, the expected value of the number of indenpendent sets in this subcase is $2E(\tilde{i}_{n - 2})$.
\vskip 5 pt

\noindent \textbf{Subscase 2.4:} The independent sets that contain  neither $h_{2}$ nor  $h_{6}$.\\
\indent The number of independent sets in this subcase equals the number of independent sets of the union of $R_{n - 1}(a, b, c)$ and the path $P_{3} = h_{3}h_{4}h_{5}$. Clearly, $P_{3}$ has $5$ independent sets. So, the expected value of the number of independent sets in this subcase is $5E(i_{n - 1})$.
\vskip 5 pt

\indent From the four subcases, we have that there are $8bE(\tilde{i}_{n - 2}) + 5bE(i_{n - 1})$.
\vskip 5 pt

\noindent \textbf{Case 3:} $h_{1}$ is at distance three from $k_{1}$.\\
\indent Thus, $h_{1}$ is $k_{4}$. This case occurs with the probability $c$. We count the number of independent sets due to the following subcases.
\vskip 5 pt

\noindent \textbf{Subscase 3.1:} The independent sets that contain $h_{2}$ and do not contain $h_{6}$.\\
\indent For any independent sets that contain $h_{2}$ and does not contain $h_{6}$, it cannot contain $h_{1}$ and $h_{3}$. Thus, the number of independent sets in this subcase equals the number of independent sets of the union of $\hat{R}_{n - 2}(a, b, c)$ and the path $P_{2} = h_{4}h_{5}$. Clearly, $P_{2}$ has $3$ independent sets. So, the expected value of the number of independent sets in this subcase is $3E(\hat{i}_{n - 2})$.
\vskip 5 pt

\noindent \textbf{Subscase 3.2:} The independent sets that contain $h_{6}$ and do not contain $h_{2}$.\\
\indent For any independent sets that contain $h_{6}$ and does not contain $h_{2}$, it cannot contain $h_{1}$ and $h_{5}$. Thus, the number of independent sets in this subcase equals the number of independent sets of the union of $\hat{R}_{n - 2}(a, b, c)$ and the path $P_{2} = h_{3}h_{4}$. Clearly, $P_{2}$ has $3$ independent sets. So, the expected value of the number of independent sets in this subcase is $3E(\hat{i}_{n - 2})$.
\vskip 5 pt

\noindent \textbf{Subscase 3.3:} The independent sets that contain both $h_{2}$ and  $h_{6}$.\\
\indent For any independent set that contains both $h_{2}$ and $h_{6}$, it cannot contain $h_{1}$,$h_{3}$ and $h_{5}$. Thus, the number of independent sets in this subcase equals the number of independent sets of the union of $\hat{R}_{n - 2}(a, b, c)$ and the vertex $h_{4}$. Clearly, the vertex $h_{4}$ has $2$ independent sets. So, the expected value of the number of indenpendent sets in this subcase is $2E(\hat{i}_{n - 2})$.
\vskip 5 pt

\noindent \textbf{Subscase 3.4:} The independent sets that contain  neither $h_{2}$ nor  $h_{6}$.\\
\indent The number of independent sets in this subcase equals the number of independent sets of the union of $R_{n - 1}(a, b, c)$ and the path $P_{3} = h_{3}h_{4}h_{5}$. Clearly, $P_{3}$ has $5$ independent sets. So, the expected value of the number of independent sets in this subcase is $5E(i_{n - 1})$.
\vskip 5 pt

\indent From the four subcases, we have that there are $8cE(\hat{i}_{n - 2}) + 5cE(i_{n - 1})$.
 Therefore, from Cases $1, 2$ and $3$, we have that
\begin{align*}
E(i_{n}) &= 8aE(i'_{n - 2}) + 5aE(i_{n - 1}) + 8bE(\tilde{i}_{n - 2}) + 5bE(i_{n - 1}) + 8cE(\hat{i}_{n - 2}) + 5cE(i_{n - 1})\notag\\
         &= 5(a + b + c)E(i_{n - 1}) + 8aE(i'_{n - 2}) + 8bE(\tilde{i}_{n - 2}) + 8cE(\hat{i}_{n - 2})\notag\\
         &= 5E(i_{n - 1}) + 8aE(i'_{n - 2}) + 8bE(\tilde{i}_{n - 2}) + 8cE(\hat{i}_{n - 2}).\notag
\end{align*}

\noindent For $n \geq 2$, we multiply $x^{n}$ throughout the above equation and sum over all $n$. We have that
\begin{align*}
\sum^{\infty}_{n = 2}E(i_{n})x^{n} = 5\sum^{\infty}_{n = 2}E(i_{n - 1})x^{n} + 8a\sum^{\infty}_{n = 2}E(i'_{n - 2})x^{n} + 8b\sum^{\infty}_{n = 2}E(\tilde{i}_{n - 2})x^{n} + 8c\sum^{\infty}_{n = 2}E(\hat{i}_{n - 2})x^{n}\notag
\end{align*}
\noindent which implies that
\begin{align*}
I(x) - E(i_{0}) - E(i_{1})x = 5x(I(x) - E(i_{0})) + 8ax^{2}I'(x) + 8bx^{2}\tilde{I}(x) + 8cx^{2}\hat(I)(x).\notag
\end{align*}
It can be checked that $E(i_{0}) = 1$ and $E(i_{1}) = 18$. Hence, 
\begin{align*}
I(x) = 1 + 13x + 5xI(x) + 8ax^{2}I'(x) + 8bx^{2}\tilde{I}(x) + 8cx^{2}\hat{I}(x).
\end{align*}
This proves Equation (\ref{eqm5}).
\qed

\noindent \emph{Proof of Equation} (\ref{eqm6}). Recall that $R'_{n}(a, b, c)$ is obtained from $R_{n}(a, b, c)$ and a path $P_{5} = x_{1}x_{2}x_{3}x_{4}x_{5}$ by identifying $x_{1}$ with a vertex of the $n^{th}$ hexagon of $R_{n}(a, b, c)$. We name the vertices of the $n^{th}$ hexagon $H_{n}$ of $R'_{n}(a, b, c)$ by $h_{1}, h_{2}, ..., h_{6}$ clockwise with $h_{1}$ is the cut vertex containing in $H_{n}$ (and in the $(n - 1)^{th}$ hexagon $H_{n - 1}$). By the definition of $R'_{n}(a, b, c)$, we have that $x_{1} \neq h_{1}$. We distinguish three cases due to the distance between $x_{1}$ and $h_{1}$. 
\vskip 5 pt

\noindent \textbf{Case 1:} $x_{1}$ is adjacent to $h_{1}$.\\
\indent Thus, $x_{1}$ is either $h_{2}$ or $h_{6}$. Without loss of generality, we let $x_{1} = h_{2}$. This case occurs with the probability $a$. We count the number of independent sets due to the following subcases.
\vskip 5 pt

\noindent \textbf{Subscase 1.1:} The independent sets that contain $x_{2}$.\\
\indent For any independent set that contains $x_{2}$, it cannot contain $x_{1}$ and $x_{3}$. Thus, the number of independent sets in this subcase equals the number of independent sets of the union of $R'_{n -1}(a, b, c)$ and the path $x_{4}x_{5}$. Clearly, $x_{4}x_{5}$ has $3$ independent sets. So, the expected value of the number of independent sets in this subcase is $3E(i'_{n - 1})$.
\vskip 5 pt

\noindent \textbf{Subscase 1.2:} The independent sets that do not contain $x_{2}$.\\
\indent The number of independent sets in this subcase equals the number of independent sets of the union of $R_{n}(a, b, c)$ and the path $x_{3}x_{4}x_{5}$. Clearly, $x_{3}x_{4}x_{5}$ has $5$ independent sets. So, the expected value of the number of independent sets in this subcase is $5E(i_{n})$.
\vskip 5 pt

\indent From the two subcases, we have that there are $3aE(i'_{n - 1}) + 5aE(i_{n})$ from Case 1.
\vskip 5 pt

\noindent \textbf{Case 2:} $x_{1}$ is at distance two from $h_{1}$.\\
\indent Thus, $x_{1}$ is either $h_{3}$ or $h_{5}$. Without loss of generality, we let $x_{1} = h_{3}$. This case occurs with the probability $b$. We count the number of independent sets due to the following subcases.
\vskip 5 pt

\noindent \textbf{Subscase 2.1:} The independent sets that contain $x_{2}$.\\
\indent For any independent set that contains $x_{2}$, it cannot contain $x_{1}$ and $x_{3}$. Thus, the number of independent sets in this subcase equals the number of independent sets of the union of $\tilde{R}_{n -1}(a, b, c)$ and the path $x_{4}x_{5}$. Clearly, $x_{4}x_{5}$ has $3$ independent sets. So, the expected value of the number of independent sets in this subcase is $3E(\tilde{i}_{n - 1})$.
\vskip 5 pt

\noindent \textbf{Subscase 2.2:} The independent sets that do not contain $x_{2}$.\\
\indent The number of independent sets in this subcase equals the number of independent sets of the union of $R_{n}(a, b, c)$ and the path $x_{3}x_{4}x_{5}$. Clearly, $x_{3}x_{4}x_{5}$ has $5$ independent sets. So, the expected value of the number of independent sets in this subcase is $5E(i_{n})$.
\vskip 5 pt

\indent From the two subcases, we have that there are $3bE(\tilde{i}_{n - 1}) + 5bE(i_{n})$ from Case 2.
\vskip 5 pt

\noindent \textbf{Case 3:} $x_{1}$ is at distance three from $h_{1}$.\\
\indent Thus, $x_{1}$ is $h_{4}$. This case occurs with the probability $c$. We count the number of independent sets due to the following subcases.
\vskip 5 pt

\noindent \textbf{Subscase 3.1:} The independent sets that contain $x_{2}$.\\
\indent For any independent set that contains $x_{2}$, it cannot contain $x_{1}$ and $x_{3}$. Thus, the number of independent sets in this subcase equals the number of independent sets of the union of $\hat{R}_{n -1}(a, b, c)$ and the path $x_{4}x_{5}$. Clearly, $x_{4}x_{5}$ has $3$ independent sets. So, the expected value of the number of independent sets in this subcase is $3E(\hat{i}_{n - 1})$.
\vskip 5 pt

\noindent \textbf{Subscase 3.2:} The independent sets that do not contain $x_{2}$.\\
\indent The number of independent sets in this subcase equals the number of independent sets of the union of $R_{n}(a, b, c)$ and the path $x_{3}x_{4}x_{5}$. Clearly, $x_{3}x_{4}x_{5}$ has $5$ independent sets. So, the expected value of the number of independent sets in this subcase is $5E(i_{n})$.
\vskip 5 pt

\indent From the two subcases, we have that there are $3cE(\hat{i}_{n - 1}) + 5cE(i_{n})$ from Case 3.
 Therefore, from Cases $1, 2$ and $3$, we have that
\begin{align*}
E(i'_{n}) &= 3aE(i'_{n -1}) + 5aE(i_{n}) + 3bE(\tilde{i}_{n - 1}) + 5bE(i_{n}) + 3cE(\hat{i}_{n - 1}) + 5cE(i_{n})\notag\\
         &= 5(a + b + c)E(i_{n}) + 3aE(i'_{n - 1}) + 3bE(\tilde{i}_{n - 1}) + 3cE(\hat{i}_{n - 1})\notag\\
         &= 5E(i_{n}) + 3aE(i'_{n - 1}) + 3bE(\tilde{i}_{n - 1}) + 3cE(\hat{i}_{n - 1}).\notag
\end{align*}
\vskip 5 pt

\noindent For $n \geq 1$, we multiply $x^{n}$ throughout the above equation and sum over all $n$. We have that
\begin{align*}
\sum^{\infty}_{n = 1}E(i'_{n})x^{n} = 5\sum^{\infty}_{n = 1}E(i_{n})x^{n} + 3a\sum^{\infty}_{n = 1}E(i'_{n - 1})x^{n} + 3b\sum^{\infty}_{n = 1}E(\tilde{i}_{n - 1})x^{n} + 3c\sum^{\infty}_{n = 1}E(\hat{i}_{n - 1})x^{n}\notag
\end{align*}
\noindent which implies that
\begin{align*}
I'(x) - E(i'_{0})  = 5(I(x) - E(i_{0})) + 3axI'(x) + 3bx\tilde{I}(x) + 3cx\hat(I)(x).\notag
\end{align*}
It can be checked that $E(i_{0}) = 1$ and $E(i'_{0}) = 13$. Hence, 
\begin{align*}
I'(x) = 8+  5I(x) + 3axI'(x) + 3bx\tilde{I}(x) + 3cx\hat{I}(x).
\end{align*}
This proves Equation (\ref{eqm6}).
\qed

\noindent \emph{Proof of Equation} (\ref{eqm7}). Recall that $\tilde{R}_{n}(a, b, c)$ is obtained from $R_{n}(a, b, c)$ and a path $P_{5} = x_{1}x_{2}x_{3}x_{4}x_{5}$ by identifying $x_{2}$ with a vertex of the $n^{th}$ hexagon of $R_{n}(a, b, c)$. We name the vertices of the $n^{th}$ hexagon $H_{n}$ of $\tilde{R}_{n}(a, b, c)$ by $h_{1}, h_{2}, ..., h_{6}$ clockwise with $h_{1}$ is the cut vertex containing in $H_{n}$ (and in the $(n - 1)^{th}$ hexagon $H_{n - 1}$). By the definition of $\tilde{R}_{n}(a, b, c)$, we have that $x_{2} \neq h_{1}$. We distinguish three cases due to the distance between $x_{2}$ and $h_{1}$. 
\vskip 5 pt

\noindent \textbf{Case 1:} $x_{2}$ is adjacent to $h_{1}$.\\
\indent Thus, $x_{2}$ is either $h_{2}$ or $h_{6}$. Without loss of generality, we let $x_{2} = h_{2}$. This case occurs with the probability $a$. We count the number of independent sets due to the following subcases.
\vskip 5 pt

\noindent \textbf{Subscase 1.1:} The independent sets that contain $x_{1}$ and do not contain $x_{3}$.\\
\indent For any independent set that contains $x_{1}$ and does not contain $x_{3}$, it cannot contain the vertex $x_{2}$. Thus,the number of independent sets in this subcase equals the number of independent sets of the union of $R'_{n -1}(a, b, c)$ and the path $x_{4}x_{5}$. Clearly, $x_{4}x_{5}$ has $3$ independent sets. So, the expected value of the number of independent sets in this subcase is $3E(i'_{n - 1})$.
\vskip 5 pt

\noindent \textbf{Subscase 1.2:} The independent sets that contain $x_{3}$ and do not contain $x_{1}$.\\
\indent For any independent set that contains $x_{3}$ and does not contain $x_{1}$, it cannot contain the vertex $x_{2}$ and $x_{4}$. Thus, the number of independent sets in this subcase equals the number of independent sets of the union of $R'_{n -1}(a, b, c)$ and the vertex $x_{5}$. Clearly, the vertex $x_{5}$ has $2$ independent sets. So, the expected value of the number of independent sets in this subcase is $2E(i'_{n - 1})$.
\vskip 5 pt

\noindent \textbf{Subscase 1.3:} The independent sets that contain both $x_{1}$ and  $x_{3}$.\\
\indent For any independent set that contains both $x_{1}$ and $x_{3}$, it cannot contain the vertex $x_{2}$ and $x_{4}$. Thus, the number of independent sets in this subcase equals the number of independent sets of the union of $R'_{n - 1}(a, b, c)$ and the vertex $h_{5}$. Clearly, the vertex $h_{5}$ has $2$ independent sets. So, the expected value of the number of indenpendent sets in this subcase is $2E(i'_{n - 1})$.
\vskip 5 pt

\noindent \textbf{Subscase 1.4:} The independent sets that contain  neither $x_{1}$ nor  $x_{3}$.\\
\indent The number of independent sets in this subcase equals the number of independent sets of the union of $R_{n }(a, b, c)$ and the path $P_{2} = x_{4}x_{5}$. Clearly, $P_{2}$ has $3$ independent sets. So, the expected value of the number of independent sets in this subcase is $3E(i_{n })$.
\vskip 5 pt

\indent From the four subcases, we have that there are $7aE(i'_{n - 1}) + 3aE(i_{n })$.
\vskip 5 pt

\noindent \textbf{Case 2:} $x_{2}$ is at distance two from $h_{1}$.\\
\indent Thus, $x_{2}$ is either $h_{3}$ or $h_{5}$. Without loss of generality, we let $x_{2} = h_{3}$. This case occurs with the probability $b$. We count the number of independent sets due to the following subcases.
\vskip 5 pt

\noindent \textbf{Subscase 2.1:} The independent sets that contain $x_{1}$ and do not contain $x_{3}$.\\
\indent For any independent set that contains $x_{1}$ and does not contain $x_{3}$, it cannot contain the vertex $x_{2}$. Thus, the number of independent sets in this subcase equals the number of independent sets of the union of $\tilde{R}_{n -1}(a, b, c)$ and the path $x_{4}x_{5}$. Clearly, $x_{4}x_{5}$ has $3$ independent sets. So, the expected value of the number of independent sets in this subcase is $3E(\tilde{i}_{n - 1})$.
\vskip 5 pt

\noindent \textbf{Subscase 2.2:} The independent sets that contain $x_{3}$ and do not contain $x_{1}$.\\
\indent For any independent set that contains $x_{3}$ and does not contain $x_{1}$, it cannot contain the vertex $x_{2}$ and $x_{4}$. Thus, the number of independent sets in this subcase equals the number of independent sets of the union of $\tilde{R}_{n -1}(a, b, c)$ and the vertex $x_{5}$. Clearly, the vertex $x_{5}$ has $2$ independent sets. So, the expected value of the number of independent sets in this subcase is $2E(\tilde{i}_{n - 1})$.
\vskip 5 pt

\noindent \textbf{Subscase 2.3:} The independent sets that contain both $x_{1}$ and $x_{3}$.\\
\indent For any independent set that contains both $x_{1}$ and  $x_{3}$, it cannot contain the vertex $x_{2}$ and $x_{4}$. Thus, the number of independent sets in this subcase equals the number of independent sets of the union of $\tilde{R}_{n -1}(a, b, c)$ and the vertex $x_{5}$. Clearly, the vertex $x_{5}$ has $2$ independent sets. So, the expected value of the number of independent sets in this subcase is $2E(\tilde{i}_{n - 1})$.
\vskip 5 pt

\noindent \textbf{Subscase 2.4:} The independent sets that contain  neither $x_{1}$ nor  $x_{3}$.\\
\indent The number of independent sets in this subcase equals the number of independent sets of the union of $R_{n}(a, b, c)$ and the path $P_{2} = x_{4}x_{5}$. Clearly, $P_{2}$ has $3$ independent sets. So, the expected value of the number of independent sets in this subcase is $3E(i_{n})$.
\vskip 5 pt

\indent From the four subcases, we have that there are $7bE(i'_{n - 1}) + 3bE(i_{n })$.
\vskip 5 pt

\noindent \textbf{Case 3:} $x_{2}$ is at distance three from $h_{1}$.\\
\indent Thus, $x_{2}$ is  $h_{4}$.  This case occurs with the probability $c$. We count the number of independent sets due to the following subcases.
\vskip 5 pt

\noindent \textbf{Subscase 3.1:} The independent sets that contain $x_{1}$ and do not contain $x_{3}$.\\
\indent For any independent set that contains $x_{1}$ and does not contain $x_{3}$, it cannot contain the vertex $x_{2}$. Thus, the number of independent sets in this subcase equals the number of independent sets of the union of $\hat{R}_{n -1}(a, b, c)$ and the path $x_{4}x_{5}$. Clearly, $x_{4}x_{5}$ has $3$ independent sets. So, the expected value of the number of independent sets in this subcase is $3E(\hat{i}_{n - 1})$.
\vskip 5 pt

\noindent \textbf{Subscase 3.2:} The independent sets that contain $x_{3}$ and do not contain $x_{1}$.\\
\indent For any independent set that contains $x_{3}$ and does not contain $x_{1}$, it cannot contain the vertex $x_{2}$ and $x_{4}$. Thus, the number of independent sets in this subcase equals the number of independent sets of the union of $\hat{R}_{n -1}(a, b, c)$ and the vertex $x_{5}$. Clearly, the vertex $x_{5}$ has $2$ independent sets. So, the expected value of the number of independent sets in this subcase is $2E(\hat{i}_{n - 1})$.
\vskip 5 pt

\noindent \textbf{Subscase 3.3:} The independent sets that contain both $x_{1}$ and $x_{3}$.\\
\indent For any independent set that contains both $x_{1}$ and  $x_{3}$, it cannot contain the vertex $x_{2}$ and $x_{4}$. Thus, the number of independent sets in this subcase equals the number of independent sets of the union of $\hat{R}_{n -1}(a, b, c)$ and the vertex $x_{5}$. Clearly, the vertex $x_{5}$ has $2$ independent sets. So, the expected value of the number of independent sets in this subcase is $2E(\hat{i}_{n - 1})$.
\vskip 5 pt

\noindent \textbf{Subscase 3.4:} The independent sets that contain  neither $x_{1}$ nor  $x_{3}$.\\
\indent The number of independent sets in this subcase equals the number of independent sets of the union of $R_{n}(a, b, c)$ and the path $P_{2} = x_{4}x_{5}$. Clearly, $P_{2}$ has $3$ independent sets. So, the expected value of the number of independent sets in this subcase is $3E(i_{n})$.
\vskip 5 pt

\indent From the four subcases, we have that there are $7cE(i'_{n - 1}) + 3cE(i_{n })$. Therefore, from Cases $1, 2$ and $3$, we have that
\begin{align*}
E(\tilde{i}_{n}) &= 7aE(i'_{n -1}) + 3aE(i_{n}) + 7bE(\tilde{i}_{n - 1}) + 3bE(i_{n}) + 7cE(\hat{i}_{n - 1}) + 3cE(i_{n})\notag\\
         &= 3(a + b + c)E(i_{n}) + 7aE(i'_{n - 1}) + 7bE(\tilde{i}_{n - 1}) + 7cE(\hat{i}_{n - 1})\notag\\
         &= 3E(i_{n}) + 7aE(i'_{n - 1}) + 7bE(\tilde{i}_{n - 1}) + 7cE(\hat{i}_{n - 1}).\notag
\end{align*}
\vskip 5 pt

\noindent For $n \geq 1$, we multiply $x^{n}$ throughout the above equation and sum over all $n$. We have that
\begin{align*}
\sum^{\infty}_{n = 1}E(\tilde{i}_{n})x^{n} = 3\sum^{\infty}_{n = 1}E(i_{n})x^{n} + 7a\sum^{\infty}_{n = 1}E(i'_{n - 1})x^{n} + 7b\sum^{\infty}_{n = 1}E(\tilde{i}_{n - 1})x^{n} + 7c\sum^{\infty}_{n = 1}E(\hat{i}_{n - 1})x^{n}\notag
\end{align*}
\noindent which implies that
\begin{align*}
\tilde{I}(x) - E(\tilde{i}_{0})  = 3(I(x) - E(i_{0})) + 7axI'(x) + 7bx\tilde{I}(x) + 7cx\hat(I)(x).\notag
\end{align*}
It can be checked that $E(i_{0}) = 1$ and $E(\tilde{i}_{0}) = 13$. Hence, 
\begin{align*}
\tilde{I}(x) = 10 + 3I(x) + 7axI'(x) + 7bx\tilde{I}(x) + 7cx\hat{I}(x).
\end{align*}
This proves Equation (\ref{eqm7}).
\qed
\vskip 5 pt

\noindent \emph{Proof of Equation} (\ref{eqm8}). Recall that $\hat{R}_{n}(a, b, c)$ is obtained from $R_{n}(a, b, c)$ and a path $P_{5} = x_{1}x_{2}x_{3}x_{4}x_{5}$ by identifying ${x_3}$ with a vertex of the $n^{th}$ hexagon of $R_{n}(a, b, c)$. We name the vertices of the $n^{th}$ hexagon $H_{n}$ of $\hat{R}_{n}(a, b, c)$ by $h_{1}, h_{2}, ..., h_{6}$ clockwise with $h_{1}$ is the cut vertex containing in $H_{n}$ (and in the $(n - 1)^{th}$ hexagon $H_{n - 1}$). By the definition of $\hat{R}_{n}(a, b, c)$, we have that $x_{3} \neq h_{1}$. We distinguish three cases due to the distance between $x_{3}$ and $h_{1}$. 
\vskip 5 pt

\noindent \textbf{Case 1:} $x_{3}$ is adjacent to $h_{1}$.\\
\indent Thus, $x_{3}$ is either $h_{2}$ or $h_{6}$. Without loss of generality, we let $x_{3} = h_{2}$. This case occurs with the probability $a$. We count the number of independent sets due to the following subcases.
\vskip 5 pt

\noindent \textbf{Subscase 1.1:} The independent sets that contain $x_{2}$ and do not contain $x_{4}$.\\
\indent For any independent set that contains $x_{2}$ and does not contain $x_{4}$, it cannot contain the vertices $x_{1}$ and $x_{3}$. Thus, the number of independent sets in this subcase equals the number of independent sets of the union of $R'_{n -1}(a, b, c)$ and the vertex $x_{5}$. Clearly, the vertex $x_{5}$ has $2$ independent sets. So, the expected value of the number of independent sets in this subcase is $2E(i'_{n - 1})$.
\vskip 5 pt

\noindent \textbf{Subscase 1.2:} The independent sets that contain $x_{4}$ and do not contain $x_{2}$.\\
\indent For any independent set that contains $x_{4}$ and does not contain $x_{2}$, it cannot contain the vertices  $x_{3}$ and $x_{5}$. Thus, the number of independent sets in this subcase equals the number of independent sets of the union of $R'_{n -1}(a, b, c)$ and the vertex $x_{1}$. Clearly, the vertex $x_{1}$ has $2$ independent sets. So, the expected value of the number of independent sets in this subcase is $2E(i'_{n - 1})$.
\vskip 5 pt

\noindent \textbf{Subscase 1.3:} The independent sets that contain both $x_{2}$ and  $x_{4}$.\\
\indent For any independent set that contains both $x_{2}$ and $x_{4}$, it cannot contain the vertex $x_{3}$, $x_{4}$ and $x_{5}$. Thus, the number of independent sets in this subcase equals the number of independent sets of $R'_{n - 1}(a, b, c)$. So, the expected value of the number of independent sets in this subcase is $E(i'_{n - 1})$.
\vskip 5 pt

\noindent \textbf{Subscase 1.4:} The independent sets that contain  neither $x_{2}$ nor  $x_{4}$.\\
\indent The number of independent sets in this subcase equals the number of independent sets of the union of $R_{n }(a, b, c)$ and the vertices  $ x_{1}$ and $x_{5}$. Clearly, the vertices  $ x_{1}$ and $x_{5}$ has $4$ independent sets. So, the expected value of the number of independent sets in this subcase is $4E(i_{n})$.
\vskip 5 pt

\indent From the four subcases, we have that there are $5aE(i'_{n - 1}) + 4aE(i_{n })$.
\vskip 5 pt

\noindent \textbf{Case 2:} $x_{3}$ is at distance two from $h_{1}$.\\
\indent Thus, $x_{3}$ is either $h_{3}$ or $h_{5}$. Without loss of generality, we let $x_{3} = h_{3}$. This case occurs with the probability $b$. We count the number of independent sets due to the following subcases.
\vskip 5 pt

\noindent \textbf{Subscase 2.1:} The independent sets that contain $x_{2}$ and do not contain $x_{4}$.\\
\indent For any independent set that contains $x_{2}$ and does not contain $x_{4}$, it cannot contain the vertex $x_{1}$ and $x_{3}$. Thus, the number of independent sets in this subcase equals the number of independent sets of the union of $\tilde{R}_{n -1}(a, b, c)$ and the vertex $x_{5}$. Clearly, the vertex $x_{5}$ has $2$ independent sets. So, the expected value of the number of independent sets in this subcase is $2E(\tilde{i}_{n - 1})$.
\vskip 5 pt

\noindent \textbf{Subscase 2.2:} The independent sets that contain $x_{4}$ and do not contain $x_{2}$.\\
\indent For any independent set that contains $x_{4}$ and does not contain $x_{2}$, it cannot contain the vertex $x_{3}$ and $x_{5}$. Thus, the number of independent sets in this subcase equals the number of independent sets of the union of $\tilde{R}_{n -1}(a, b, c)$ and the vertex $x_{1}$. Clearly, the vertex $x_{1}$ has $2$ independent sets. So, the expected value of the number of independent sets in this subcase is $2E(\tilde{i}_{n - 1})$.
\vskip 5 pt

\noindent \textbf{Subscase 2.3:} The independent sets that contain both $x_{2}$ and  $x_{4}$.\\
\indent For any independent set that contains both $x_{2}$ and $x_{4}$, it cannot contain the vertex $x_{3}$, $x_{4}$ and $x_{5}$. Thus,the number of independent sets in this subcase equals the number of independent sets of  $\tilde{R}_{n - 1}(a, b, c)$. So, the expected value of the number of independent sets in this subcase is $E(\tilde{i}_{n - 1})$.
\vskip 5 pt

\noindent \textbf{Subscase 2.4:} The independent sets that contain  neither $x_{2}$ nor  $x_{4}$.\\
\indent The number of independent sets in this subcase equals the number of independent sets of the union of $R_{n }(a, b, c)$ and the vertices  $ x_{1}$ and $x_{5}$. Clearly, the vertices  $ x_{1}$ and $x_{5}$ has $4$ independent sets. So, the expected value of the number of independent sets in this subcase is $4E(i_{n})$.
\vskip 5 pt

\indent From the four subcases, we have that there are $5bE(i'_{n - 1}) + 4bE(i_{n })$.
\vskip 5 pt

\noindent \textbf{Case 3:} $x_{3}$ is at distance three from $h_{1}$.\\
\indent Thus, $x_{3}$ is  $h_{4}$.  This case occurs with the probability $c$. We count the number of independent sets due to the following subcases.
\vskip 5 pt

\noindent \textbf{Subscase 3.1:} The independent sets that contain $x_{2}$ and do not contain $x_{4}$.\\
\indent For any independent set that contains $x_{2}$ and does not contain $x_{4}$, it cannot contain  $x_{1}$ and $x_{3}$. Thus, the number of independent sets in this subcase equals the number of independent sets of the union of $\hat{R}_{n -1}(a, b, c)$ and the vertex $x_{5}$. Clearly, the vertex $x_{5}$ has $2$ independent sets. So, the expected value of the number of independent sets in this subcase is $2E(\hat{i}_{n - 1})$.
\vskip 5 pt

\noindent \textbf{Subscase 3.2:} The independent sets that contain $x_{4}$ and do not contain $x_{2}$.\\
\indent For any independent set that contains $x_{4}$ and does not contain $x_{2}$, it cannot contain the vertex $x_{3}$ and $x_{5}$. Thus, the number of independent sets in this subcase equals the number of independent sets of the union of $\hat{R}_{n -1}(a, b, c)$ and the vertex $x_{1}$. Clearly, the vertex $x_{1}$ has $2$ independent sets. So, the expected value of the number of independent sets in this subcase is $2E(\hat{i}_{n - 1})$.
\vskip 5 pt

\noindent \textbf{Subscase 3.3:} The independent sets that contain both $x_{2}$ and  $x_{4}$.\\
\indent For any independent set that contains both $x_{2}$ and $x_{4}$, it cannot contain the vertex $x_{3}$, $x_{4}$ and $x_{5}$. Thus, the number of independent sets in this subcase equals the number of independent sets of  $\hat{R}_{n - 1}(a, b, c)$. So, the expected value of the number of independent sets in this subcase is $E(\hat{i}_{n - 1})$.
\vskip 5 pt

\noindent \textbf{Subscase 3.4:} The independent sets that contain  neither $x_{2}$ nor  $x_{4}$.\\
\indent The number of independent sets in this subcase equals the number of independent sets of the union of $R_{n }(a, b, c)$ and the vertices  $ x_{1}$ and $x_{5}$. Clearly, the vertices  $ x_{1}$ and $x_{5}$ has $4$ independent sets. So, the expected value of the number of independent sets in this subcase is $4E(i_{n})$.
\vskip 5 pt

\indent From the four subcases, we have that there are $5cE(i'_{n - 1}) + 4cE(i_{n })$. Therefore, from Cases $1, 2$ and $3$, we have that
\begin{align*}
E(\hat{i}_{n}) &= 5aE(i'_{n -1}) + 4aE(i_{n}) + 5bE(\tilde{i}_{n - 1}) + 4bE(i_{n}) + 5cE(\hat{i}_{n - 1}) + 4cE(i_{n})\notag\\
         &= 4(a + b + c)E(i_{n}) + 5aE(i'_{n - 1}) + 5bE(\tilde{i}_{n - 1}) + 5cE(\hat{i}_{n - 1})\notag\\
         &= 4E(i_{n}) + 5aE(i'_{n - 1}) + 5bE(\tilde{i}_{n - 1}) + 5cE(\hat{i}_{n - 1}).\notag
\end{align*}
\vskip 5 pt

\noindent For $n \geq 1$, we multiply $x^{n}$ throughout the above equation and sum over all $n$. We have that
\begin{align*}
\sum^{\infty}_{n = 1}E(\hat{i}_{n})x^{n} = 4\sum^{\infty}_{n = 1}E(i_{n})x^{n} + 5a\sum^{\infty}_{n = 1}E(i'_{n - 1})x^{n} + 5b\sum^{\infty}_{n = 1}E(\tilde{i}_{n - 1})x^{n} + 5c\sum^{\infty}_{n = 1}E(\hat{i}_{n - 1})x^{n}\notag
\end{align*}
\noindent which implies that
\begin{align*}
\hat{I}(x) - E(\hat{i}_{0})  = 4(I(x) - E(i_{0})) + 5axI'(x) + 5bx\tilde{I}(x) + 5cx\hat(I)(x).\notag
\end{align*}
It can be checked that $E(i_{0}) = 1$ and $E(\hat{i}_{0}) = 13$. Hence, 
\begin{align*}
\hat{I}(x) = 9 + 4I(x) + 5axI'(x) + 5bx\tilde{I}(x) + 5cx\hat{I}(x).
\end{align*}
This proves Equation (\ref{eqm8}).
\qed
\vskip 5 pt

\indent By Equations (\ref{eqm5}), (\ref{eqm6}), (\ref{eqm7}) and (\ref{eqm8}), it can be solved that 
\begin{align*}
I(x) =  \frac{1 + 13x - 3ax - 7bx - 5cx + 25ax^{2} - 11bx^{2} + 7cx^{2}}{1 - 5x - 3ax - 7bx - 5cx - 25ax^{2} + 11bx^{2} - 7cx^{2}}.
\end{align*}
This proves Theorem \ref{theorem2}.


\noindent \textbf{Acknowledgment:} The first author has been supported by Petchra Pra Jom Klao Master Scholarship from King Mongkut's University of Technology Thonburi~(21/2565).

\

\end{document}